\documentclass[5p,times]{elsarticle}
\journal{Journal of Process Control}

\usepackage{amsfonts}
\usepackage{float}
\usepackage{amsmath}
\usepackage{graphicx}
\usepackage{amsthm}
\usepackage{epstopdf}
\usepackage{subfigure}
\usepackage{color}
\usepackage{epsfig}
\usepackage{multirow}
\usepackage{fancyhdr}
\usepackage{psfrag}
\usepackage{dsfont}
\usepackage{amssymb}
\usepackage{xcolor}
\usepackage[normalem]{ulem}

\usepackage{moreverb}
\usepackage{graphicx}
\usepackage{float}
\usepackage{hyperref}
\usepackage{tabularx}
\usepackage{array}
\usepackage{amssymb}
\usepackage{amsmath}

\newtheorem{theorem}{Theorem}

\newtheorem{lemma}[theorem]{Lemma}

\newtheorem{remark}[theorem]{Remark}
\newtheorem{proposition}[theorem]{Proposition}

\usepackage[figuresright]{rotating}
\begin{document}

\begin{frontmatter}

\title{Boundary-controlled irreversible port--Hamiltonian modeling and
structure-preserving discretization of tubular reactors}
\tnotetext[mytitlenote]{}

\author[First]{Weijun. Zhou}

\address[First]{School of Information and Electrical Engineering, Hangzhou City University, 310027, P. R. China,(e-mail: zhouwj@hzcu.edu.cn)}

\begin{abstract}

This paper presents a boundary-controlled irreversible port--Hamiltonian
model and a structure-preserving spatial discretization for a
non-isothermal tubular reactor with axial convection, variable-density
momentum dynamics, chemical reaction, mass diffusion, heat conduction,
and jacket heat transfer. Reaction, diffusion, and heat conduction are
described by irreversible port--Hamiltonian channels, whereas convection
and momentum are represented by a weighted Eulerian interconnection.
The boundary ports recover the axial energy flux and the corresponding
entropy flux. A partitioned finite-element method is then formulated in conservative
coordinates using a power-preserving change of variables. The
semi-discrete model satisfies the Hamiltonian balance and a discrete
entropy balance with non-negative entropy production due to reaction,
diffusion, heat conduction, and jacket heat transfer. Numerical results
for an open tubular reactor are validated against independently computed
continuous steady boundary-value solutions and show the transient
temperature, composition, and velocity profiles in the very-low-Mach
regime. The same discretization is also tested with reciprocal
Soret--Dufour coupling, for which the discrete entropy production remains
non-negative.

\end{abstract}


\end{frontmatter}


\section{Introduction}

Port--Hamiltonian systems (PHS) provide a unified framework for the
modeling and control of multiphysics systems by combining energy-based
representations with structure-preserving interconnections. Owing to their
intrinsic compatibility with conservation laws and power exchanges, PHS
have been successfully applied to both lumped- and distributed-parameter
systems, where the underlying dynamics are described through Dirac or
Stokes--Dirac structures
\cite{maschke1992,vanderschaft2002,legorrec2005,maschke2013,Rashad2020}.
Many engineering processes, however, involve not only reversible energy
exchange but also irreversible phenomena such as diffusion, heat conduction
and chemical reaction. Their physically meaningful description therefore
requires the simultaneous satisfaction of the first and second laws of
thermodynamics, including an explicit representation of entropy production.
Irreversible port--Hamiltonian systems (IPHS) address this requirement by
incorporating thermodynamic irreversibility while retaining the underlying
power-based interconnection structure
\cite{ramirez2012,Ramirez2013,ramirez2022overview,ramirez2022bciphs}.

For distributed irreversible systems, boundary interactions are particularly
important because they determine how energy and entropy are exchanged with
the environment and how distributed subsystems may subsequently be
interconnected. This issue motivated the development of boundary-controlled
IPHS formulations, in which the Green--Stokes identity is used to construct
power-preserving boundary variables compatible with irreversible
thermodynamics \cite{ramirez2022bciphs}. Tubular reactors provide a demanding
application of this framework: multicomponent convection, dynamic momentum,
chemical reaction, mass diffusion and heat transport evolve simultaneously,
while energy and material are exchanged through the inlet, outlet and
distributed reactor--jacket interface. The main difficulty is to combine variable-density convection with
irreversible reaction and transport while retaining both the physical
boundary-power balance and entropy-production properties after spatial
discretization.

Several developments provide important ingredients for such a formulation.
An entropy-based distributed port--Hamiltonian representation of irreversible
processes was introduced in \cite{zhou2017distributed}, where the Gibbs
relation was combined with distributed balance laws to describe irreversible
transport and reaction mechanisms. A thermodynamically consistent tubular
reactor model was subsequently developed in
\cite{zhou2021thermodynamic}, including species, energy, and entropy balances
under the local-equilibrium assumption. Related port--Hamiltonian formulations have also been developed for reactive
fluid systems. In particular, Altmann and Schulze
\cite{AltmannSchulze2017} derived a one-dimensional port--Hamiltonian
formulation of the Navier--Stokes equations for reactive flows with non-zero
boundary energy exchange. More recently, irreversible port--Hamiltonian
formulations for one-dimensional fluid systems have provided an Eulerian
framework incorporating distributed convection and momentum dynamics
\cite{Ouardi2026}. Thermodynamically consistent interconnections of
irreversible port--Hamiltonian systems have also been investigated in
\cite{RamirezLeGorrec2024}. On the numerical side, structure-preserving spatial discretization methods
have been developed for distributed port--Hamiltonian systems
\cite{Kotyczka2016,KotyczkaMaschkeLefevre2018,Rashad2020}.
In particular, the partitioned finite-element method (PFEM)
\cite{CardosoRibeiro2018,CardosoRibeiro2021,Serhani2019}
provides a weak-form construction that preserves the boundary power balance
at the semi-discrete level. To the authors' knowledge, these developments have not yet been combined
into a boundary-controlled irreversible port--Hamiltonian formulation of an
open tubular reactor that simultaneously incorporates dynamic Eulerian
convection and momentum, reaction, diffusion, and heat transport, together
with a thermodynamically consistent structure-preserving spatial
discretization.

The present work combines an irreversible port--Hamiltonian description
of reaction, mass diffusion, and heat conduction with the weighted
Eulerian interconnection for variable-density convection. Power-preserving boundary ports are constructed so that the physical
axial energy flux and entropy balance are recovered explicitly.
A conservative-coordinate PFEM discretization is then derived and shown
to preserve the semi-discrete Hamiltonian balance while retaining
non-negative entropy-production contributions. The formulation is
assessed on a dimensional open tubular reactor, with a separate
Soret--Dufour case used only as a constitutive comparison.

\section{Thermodynamic Foundations}
\label{sec:thermodynamics}

This section collects the thermodynamic relations required for the
port--Hamiltonian formulation developed in Section~\ref{sec:standard_iphs}.
Starting from the tubular-reactor model of
\cite{zhou2021thermodynamic}, the formulation is extended to variable
density and distributed momentum dynamics.

We consider a one-dimensional tubular reactor of length $L$, with spatial coordinate $z\in[0,L]$, undergoing the binary reaction $A \rightleftharpoons B$. The signed reaction rate $r$ is taken positive in the
$A\rightarrow B$ direction.

The mixture is assumed to satisfy local thermodynamic equilibrium. The
distributed physical variables are $\rho_A(z,t)$, $\rho_B(z,t)$,
$v(z,t)$, and $s(z,t)$, where $\rho_A$ and $\rho_B$ are the partial
mass densities, $v$ is the axial velocity, and $s$ is the specific
entropy. The total density and mass fractions are
\[
    \rho=\rho_A+\rho_B,
    \alpha_A=\frac{\rho_A}{\rho},
    \alpha_B=\frac{\rho_B}{\rho},
    \alpha_A+\alpha_B=1 .
\]

In contrast with the constant-density setting of
\cite{zhou2021thermodynamic}, $\rho$ is not prescribed. Indeed, if the
density were spatially and temporally constant, the total-mass balance
would imply $\partial_zv=0$, preventing $v$ from acting as a genuinely
distributed Eulerian state. Allowing $\rho$ to evolve provides the
kinematic setting required for the convective port--Hamiltonian
formulation.

The species balances are
\begin{equation}
\begin{aligned}
    \partial_t\rho_A+\partial_z(\rho_Av+J)&=-r,\\
    \partial_t\rho_B+\partial_z(\rho_Bv-J)&=r,
\end{aligned}
\label{eq:species_balances}
\end{equation}
where $J$ is the diffusive mass flux of species $A$ relative to the
barycentric motion; the corresponding diffusive flux of species $B$
is $-J$. Their sum gives
\[
    \partial_t\rho+\partial_z(\rho v)=0 .
\]

Neglecting body forces and viscous stresses, the axial momentum balance is
\begin{equation}
    \partial_t v
    +
    v\partial_zv
    +
    \frac{1}{\rho}\partial_zp
    =
    0 ,
    \label{eq:momentum}
\end{equation}
where $p$ denotes the thermodynamic pressure.

The reactor exchanges heat with a surrounding jacket at temperature
$T_j(z,t)$. No mass transfer through the reactor wall is considered.
The distributed heat input to the reacting mixture is
\[
    q_j
    =
    \kappa_j(T_j-T),
    \qquad
    \kappa_j\geq0,
\]
where $q_j>0$ denotes heat supplied to the reactor.

The entropy density is $\eta:=\rho s$, and
$\varphi:=1/\rho$ denotes the specific volume. For a binary mixture in
local thermodynamic equilibrium, the specific Gibbs relation is
\[
    du
    =
    T\,ds
    -
    p\,d\varphi
    +
    \mu_A\,d\alpha_A
    +
    \mu_B\,d\alpha_B ,
\]
where $u$ is the specific internal energy and $\mu_A,\mu_B$ are the
mass-specific chemical potentials. Since
$d\alpha_B=-d\alpha_A$, introducing the specific reaction affinity
$A_r:=\mu_A-\mu_B$ gives
\begin{equation}
    du
    =
    T\,ds
    -
    p\,d\varphi
    +
    A_r\,d\alpha_A .
    \label{eq:Gibbs_reduced}
\end{equation}

The corresponding Euler relation is
\[
    u+\frac{p}{\rho}
    =
    Ts+\alpha_A\mu_A+\alpha_B\mu_B .
\]
Hence, with the specific enthalpy
$h:=u+p/\rho$,
\[
    h
    =
    Ts+\alpha_A\mu_A+\alpha_B\mu_B .
\]
Differentiation of the Euler relation together with the Gibbs relation
gives the Gibbs--Duhem identity
\[
    dp
    =
    \rho s\,dT
    +
    \rho_A\,d\mu_A
    +
    \rho_B\,d\mu_B ,
\]
or equivalently
\[
    \frac{1}{\rho}dp
    =
    s\,dT
    +
    \alpha_A\,d\mu_A
    +
    \alpha_B\,d\mu_B .
\]

The constitutive relation between the reaction affinity and the signed
reaction rate is assumed to admit the regular factorization
\begin{equation}
    r
    =
    \Phi(A_r)
    =
    A_r\,\ell_r(A_r),
    \qquad
    \ell_r(A_r)\geq0 ,
    \label{eq:reaction_dissipation}
\end{equation}
where $\Phi(0)=0$ and $\ell_r$ admits a finite non-negative extension at
$A_r=0$. Consequently,
\[
    \sigma_r
    :=
    \frac{A_rr}{T}
    =
    \frac{A_r^2\ell_r(A_r)}{T}
    \geq0 .
\]
This regular factorization also provides the modulation function used
for the reaction channel in
Section~\ref{sec:standard_interconnections}.

Throughout this work, $q$ denotes the non-convective thermal energy flux
after the diffusive chemical-energy transport has been separated
explicitly. Hence the total non-convective energy flux is $A_rJ+q$, whereas the corresponding non-convective entropy flux is $q/T$. This
convention is used consistently in the entropy and boundary-power
balances.

The diffusion--heat contribution to the local entropy production is
\[
    -\frac{J}{T}\partial_zA_r
    +
    q\,\partial_z\!\left(\frac{1}{T}\right),
\]
with entropy-conjugate thermodynamic forces
\[
    X_d
    :=
    -\frac{1}{T}\partial_zA_r,
    \qquad
    X_q
    :=
    \partial_z\!\left(\frac{1}{T}\right).
\]

The diagonal transport case used in
Sections~\ref{sec:standard_iphs}--\ref{sec:discretization} is defined by
\begin{equation}
    J=-M_d\,\partial_zA_r,
    \qquad
    q=-k\,\partial_zT,
    \qquad
    M_d\geq0,\quad k\geq0 .
    \label{eq:diagonal_transport_laws}
\end{equation}
Equivalently, $J=M_dT X_d$ and $q=kT^2X_q$. The associated entropy
production is
\[
    \sigma_d
    =
    \frac{M_d}{T}(\partial_zA_r)^2
    \geq0,
    \qquad
    \sigma_q
    =
    \frac{k}{T^2}(\partial_zT)^2
    \geq0,
\]
and therefore $\sigma_d+\sigma_q\geq0$.

Reciprocal coupling between mass and heat transport is a standard
consequence of linear non-equilibrium thermodynamics and gives rise
to the Soret and Dufour cross-effects
\cite{degroot1962,vanbrunt2022}.

For the numerical assessment in Section~\ref{sec:numerical_thermodynamics},
a Soret--Dufour extension is considered through
\begin{equation}
    \begin{bmatrix}
        J\\ q
    \end{bmatrix}
    =
    \begin{bmatrix}
        M_dT & L_c\\
        L_c & kT^2
    \end{bmatrix}
    \begin{bmatrix}
        X_d\\ X_q
    \end{bmatrix}.
    \label{eq:onsager_matrix}
\end{equation}
Thermodynamic admissibility requires
\[
    M_d\geq0,\qquad
    k\geq0,\qquad
    L_c^2\leq M_dkT^3,
\]
which ensures
\[
    \sigma_{dq}
    =
    \frac{M_d}{T}(\partial_zA_r)^2
    +
    \frac{2L_c}{T^3}\partial_zA_r\,\partial_zT
    +
    \frac{k}{T^2}(\partial_zT)^2
    \geq0 .
\]
The case $L_c=0$ recovers the diagonal transport model used in the
IPHS formulation of Sections~\ref{sec:standard_iphs}--\ref{sec:discretization};
$L_c\neq0$ is considered only as a thermodynamically admissible
constitutive extension in the numerical comparison.

The local entropy balance follows from the Gibbs relation together with
the species, total-mass, momentum, and internal-energy balances. Adapting
the Gibbs-based derivation of \cite{zhou2021thermodynamic} to the present
variable-density model gives
\begin{equation}
    \partial_t(\rho s)
    +
    \partial_z(\rho sv)
    =
    \frac{A_rr}{T}
    -
    \frac{J}{T}\partial_zA_r
    -
    \frac{1}{T}\partial_zq
    +
    \frac{q_j}{T}.
    \label{eq:entropy_nonconservative}
\end{equation}
Using $    -\frac{1}{T}\partial_zq
    =
    -\partial_z\!\left(\frac{q}{T}\right)
    +
    q\,\partial_z\!\left(\frac{1}{T}\right)$, one obtains for the standard diagonal transport model
\begin{equation}
    \partial_t(\rho s)
    +
    \partial_z
    \left(
        \rho sv+\frac{q}{T}
    \right)
    =
    \sigma_r
    +
    \sigma_d
    +
    \sigma_q
    +
    \frac{q_j}{T}.
    \label{eq:entropy_reactor_side}
\end{equation}

The reactor-side entropy inflow associated with the jacket can be separated
into entropy supplied by the surroundings and entropy generated by the
finite temperature difference:
\begin{equation}
    \frac{q_j}{T}
    =
    \frac{q_j}{T_j}
    +
    \sigma_j ,
    \label{eq:jacket_entropy_split}
\end{equation}
where
\begin{equation}
    \sigma_j
    :=
    q_j
    \left(
        \frac{1}{T}
        -
        \frac{1}{T_j}
    \right)
    =
    \kappa_j
    \frac{(T_j-T)^2}{TT_j}
    \geq0 .
    \label{eq:jacket_entropy_production}
\end{equation}

Consequently, the complete local entropy balance of the standard model is
\begin{equation}
    \partial_t(\rho s)
    +
    \partial_z
    \left(
        \rho sv+\frac{q}{T}
    \right)
    =
    \frac{q_j}{T_j}
    +
    \sigma_r
    +
    \sigma_d
    +
    \sigma_q
    +
    \sigma_j .
    \label{eq:entropy_balance_complete}
\end{equation}
Under
\begin{equation}
    \ell_r\geq0,
    \qquad
    M_d\geq0,
    \qquad
    k\geq0,
    \qquad
    \kappa_j\geq0,
\end{equation}
the total local irreversible entropy production satisfies
\begin{equation}
    \sigma_r
    +
    \sigma_d
    +
    \sigma_q
    +
    \sigma_j
    \geq0 .
    \label{eq:entropy_inequality_complete}
\end{equation}

Equation~\eqref{eq:entropy_balance_complete} separates axial entropy
transport, jacket entropy supply, and internal irreversible entropy
production.

\section{Irreversible port--Hamiltonian formulation}
\label{sec:standard_iphs}

The objective of this section is to formulate the reacting tubular
reactor as an irreversible port--Hamiltonian system by combining the
IPHS representation of irreversible thermodynamic mechanisms
developed by Ramirez et al.~\cite{ramirez2022bciphs} with the
weighted Eulerian interconnection introduced for convective systems
in~\cite{Ouardi2026}. The reversible convective interconnection and
the irreversible reaction, diffusion, and heat-conduction channels
are distinguished explicitly below.

The structural construction concerns the diagonal transport case with $L_c=0,
    J=-M_d\,\partial_z A_r,
    q=-k\,\partial_z T$.

\subsection{Eulerian thermodynamic coordinates}
\label{sec:standard_coordinates}

Introduce the reciprocal density and the mass fraction of species $A$,
\[
    \varphi:=\frac{1}{\rho},
    \qquad
    \alpha:=\alpha_A,
    \qquad
    \alpha_B=1-\alpha .
\]
This is a change of thermodynamic state coordinates only; the spatial
coordinate $z\in[0,L]$ remains Eulerian.

With these variables, the Gibbs relation becomes
\[
    du
    =
    T\,ds
    -
    p\,d\varphi
    +
    A_r\,d\alpha ,
\]
and the Hamiltonian is
\begin{equation}
    H
    =
    \int_0^L
    \frac{1}{\varphi}
    \left(
        u(s,\varphi,\alpha)
        +
        \frac12v^2
    \right)dz .
    \label{eq:standard_H}
\end{equation}

Defining $E_c:=h+\frac12v^2$, $h=u+p\varphi$, the Hamiltonian efforts are
\begin{equation}
    e_\varphi
    =
    -\frac{E_c}{\varphi^2},
    \qquad
    e_\alpha
    =
    \frac{A_r}{\varphi},
    \qquad
    e_v
    =
    \frac{v}{\varphi},
    \qquad
    e_s
    =
    \frac{T}{\varphi}.
    \label{eq:standard_efforts}
\end{equation}

Let $x:=
    \begin{bmatrix}
        \varphi&\alpha&v
    \end{bmatrix}^{\top}$, $e_x:=
    \begin{bmatrix}
        e_\varphi&e_\alpha&e_v
    \end{bmatrix}^{\top}$.

The total entropy is
\begin{equation}
    S
    =
    \int_0^L\frac{s}{\varphi}\,dz,
    \label{eq:standard_entropy_functional}
\end{equation}
so that
\[
    \frac{\delta S}{\delta\varphi}
    =
    -\frac{s}{\varphi^2},
    \qquad
    \frac{\delta S}{\delta\alpha}=0,
    \qquad
    \frac{\delta S}{\delta v}=0,
    \qquad
    \frac{\delta S}{\delta s}
    =
    \frac1\varphi .
\]

\subsection{Convective and irreversible interconnections}
\label{sec:standard_interconnections}

\begin{proposition}[Exact convective energy and entropy transport]
\label{prop:convective_flux_exactness}
Consider the reversible convective subsystem
\[
\begin{bmatrix}
    \partial_t x\\[1mm]
    \partial_t s
\end{bmatrix}_{\!c}
=
\mathcal J_c
\begin{bmatrix}
    e_x\\[1mm]
    e_s
\end{bmatrix}.
\]
Define the local Hamiltonian and entropy densities by
\[
    \mathcal H
    =
    \frac{1}{\varphi}
    \left(
        u+\frac12v^2
    \right),
    \qquad
    \mathcal S
    =
    \frac{s}{\varphi}
    =
    \rho s .
\]
Then
\begin{equation}
    \partial_t\mathcal H
    +
    \partial_z(\rho vE_c)
    =
    0,
    \label{eq:convective_energy_exact}
\end{equation}
and
\begin{equation}
    \partial_t\mathcal S
    +
    \partial_z(\rho sv)
    =
    0.
    \label{eq:convective_entropy_exact}
\end{equation}
Consequently,
\begin{equation}
    \dot H_c
    =
    -[\rho vE_c]_0^L,
    \qquad
    \dot S_c
    =
    -[\rho sv]_0^L .
    \label{eq:convective_global_exact}
\end{equation}
Thus reversible convection generates neither internal Hamiltonian
dissipation nor internal entropy production.
\end{proposition}

\begin{proof}
Since
\[
    \partial_t\mathcal H
    =
    e_x^\top\partial_t x
    +
    e_s\partial_t s,
\]
the zero-order off-diagonal terms of $\mathcal J_c$ cancel
pointwise. The differential contribution is
\[
    e_x^\top P_1\partial_z^\varphi e_x
    =
    e_\varphi\varphi\partial_z(\varphi e_v)
    +
    e_v\varphi\partial_z(\varphi e_\varphi).
\]
Using
\[
    \varphi e_v=v,
    \qquad
    \varphi e_\varphi
    =
    -\frac{E_c}{\varphi}
    =
    -\rho E_c,
\]
gives
\[
    e_x^\top P_1\partial_z^\varphi e_x
    =
    -\rho E_c\,\partial_zv
    -
    v\,\partial_z(\rho E_c)
    =
    -\partial_z(\rho vE_c),
\]
which proves \eqref{eq:convective_energy_exact}.

For the entropy density, the first and last rows of the convective
subsystem give
\[
    \partial_t\varphi
    =
    \varphi\partial_zv
    -
    v\partial_z\varphi,
    \qquad
    \partial_t s
    =
    -C^\top G_c^\top e_x
    =
    -v\partial_zs .
\]
Therefore,
\[
    \partial_t\!\left(\frac{s}{\varphi}\right)
    =
    \frac{\partial_t s}{\varphi}
    -
    \frac{s}{\varphi^2}\partial_t\varphi
    =
    -\partial_z\!\left(\frac{sv}{\varphi}\right)
    =
    -\partial_z(\rho sv),
\]
which proves \eqref{eq:convective_entropy_exact}. Integration over
$[0,L]$ yields \eqref{eq:convective_global_exact}.
\end{proof}

The chemical reaction is represented by the active zero-order
irreversible direction
\[
    g_r=
    \begin{bmatrix}
        0 & -1 &0
    \end{bmatrix}^T,
    \qquad
    \{S|g_r|H\}
    =
    \frac{A_r}{\varphi^2}.
\]
Using the reaction factorization
\eqref{eq:reaction_dissipation}, let
\[
    \gamma_r
    :=
    \frac{\varphi^4}{T}\ell_r(A_r)
    \geq0 .
\]
Then
\begin{equation}
    R_r
    :=
    \gamma_r\{S|g_r|H\}
    =
    \frac{\varphi^2r}{T},
    \label{eq:standard_R0}
\end{equation}
and the corresponding pointwise irreversible interconnection is
\begin{equation}
    \mathcal J_r
    :=
    \begin{bmatrix}
        0 & g_rR_r\\
        -R_rg_r^\top & 0
    \end{bmatrix}.
    \label{eq:reaction_iphs_operator}
\end{equation}
Its contribution to the entropy equation is
$\varphi A_rr/T$, corresponding to the non-negative
entropy-production density $\sigma_r=A_rr/T$.

Mass diffusion is represented by the active first-order direction
\[
    g_d=
    \begin{bmatrix}
        0 & -1 & 0
    \end{bmatrix}^T,
    \qquad
    \left\{
        S
        \left|
        g_d\partial_z^\varphi
        \right|
        H
    \right\}
    =
    -\partial_zA_r .
\]
With $\gamma_d:=M_d/T\geq0$, define
\begin{equation}
    R_d
    :=
    \gamma_d
    \left\{
        S
        \left|
        g_d\partial_z^\varphi
        \right|
        H
    \right\}
    =
    -\frac{M_d}{T}\partial_zA_r
    =
    \frac JT .
    \label{eq:standard_R1}
\end{equation}
The corresponding irreversible differential interconnection is
\begin{equation}
    \mathcal J_d
    :=
    \begin{bmatrix}
        0 &
        g_d\partial_z^\varphi(R_d\,\cdot)
        \\[1mm]
        -\bigl[
            g_d\partial_z^\varphi(R_d\,\cdot)
        \bigr]^*
        &
        0
    \end{bmatrix}.
    \label{eq:diffusion_iphs_operator}
\end{equation}
Its entropy-production density is
\[
    \sigma_d
    =
    -\frac JT\,\partial_zA_r
    =
    \frac{M_d}{T}(\partial_zA_r)^2
    \geq0 .
\]

For heat conduction, retain the scalar entropy direction
\[
    g_s=\varphi,
    \qquad
    \{S|g_s|H\}
    =
    \frac1\varphi\,\partial_zT .
\]
With $\gamma_s:=\frac{k}{T^2}\geq0$, $r_s
    :=
    \gamma_s\{S|g_s|H\}
    =
    \frac{k}{\varphi T^2}\partial_zT$, the heat-conduction interconnection is
\begin{equation}
    \mathcal J_s
    :=
    \begin{bmatrix}
        0&0\\
        0&
        r_s\partial_z^\varphi(g_s\,\cdot)
        +
        g_s\partial_z^\varphi(r_s\,\cdot)
    \end{bmatrix}.
    \label{eq:standard_Js}
\end{equation}
Its entropy-production density is
\[
    \sigma_q
    =
    \frac{k}{T^2}(\partial_zT)^2
    \geq0 .
\]

Combining the reversible convective interconnection with the three
active irreversible mechanisms gives
\begin{equation}
    \begin{bmatrix}
        \partial_tx\\[1mm]
        \partial_ts
    \end{bmatrix}
    =
    \left(
        \mathcal J_c
        +
        \mathcal J_r
        +
        \mathcal J_d
        +
        \mathcal J_s
    \right)
    \begin{bmatrix}
        e_x\\[1mm]
        e_s
    \end{bmatrix}
    +
    G_jf_j .
    \label{eq:standard_complete_iphs}
\end{equation}
The internal operators are formally skew-adjoint modulo their
boundary terms. Irreversibility is encoded through the
co-energy-dependent constitutive modulations $R_r$, $R_d$, and $r_s$,
which generate the entropy-production terms associated with chemical
reaction, mass diffusion, and heat conduction.

The distributed jacket interaction is represented by
\[
    G_j
    :=
    \begin{bmatrix}
        0&0&0&\varphi
    \end{bmatrix}^T,
    \qquad
    f_j:=\frac{q_j}{T},
    \qquad
    e_j=T,
    \qquad
    e_jf_j=q_j .
\]

\begin{proposition}[Reactor balance equations]
\label{prop:standard_recovery}

Assume that the state remains in the thermodynamically admissible
domain, that $M_d\geq0$ and $k\geq0$, and that the reaction law
satisfies \eqref{eq:reaction_dissipation}. Then
\eqref{eq:standard_complete_iphs} is equivalent to
\begin{align}
    \partial_t\varphi
    &=
    \varphi\partial_zv
    -
    v\partial_z\varphi,
    \nonumber\\
    \partial_t\alpha
    &=
    -v\partial_z\alpha
    -
    \varphi r
    -
    \varphi\partial_zJ,
    \nonumber\\
    \partial_tv
    &=
    -v\partial_zv
    -
    \varphi\partial_zp,
    \label{eq:standard_recovered_system}
    \\
    \partial_ts
    &=
    -v\partial_zs
    +
    \varphi\frac{A_rr}{T}
    -
    \varphi\frac JT\partial_zA_r
    \nonumber\\
    &\quad
    -
    \varphi\partial_z\left(\frac qT\right)
    +
    \varphi\frac{k}{T^2}(\partial_zT)^2
    +
    \varphi\frac{q_j}{T}.
    \nonumber
\end{align}
Under $\rho=\frac1\varphi$, $\rho_A=\rho\alpha$, $\rho_B=\rho(1-\alpha)$, these equations recover the species, total-mass, momentum, and entropy
balances of Section~\ref{sec:thermodynamics} for $L_c=0$.
\end{proposition}

\begin{proof}
The mass and species equations follow from
$\mathcal J_c$, $\mathcal J_r$, and $\mathcal J_d$ together with
\eqref{eq:standard_R0} and \eqref{eq:standard_R1}.
For the velocity equation, the Gibbs relation gives
\begin{equation}
    dE_c
    =
    v\,dv
    +
    T\,ds
    +
    \varphi\,dp
    +
    A_r\,d\alpha,
\end{equation}
which reduces the third component of
\eqref{eq:standard_complete_iphs} to the momentum balance.

Proposition~\ref{prop:convective_flux_exactness} gives the convective entropy
transport, while the reaction, diffusion, and heat-conduction
interconnections recover the corresponding non-negative entropy-production
terms. The jacket port contributes $\varphi q_j/T$. Transforming back to
$(\rho_A,\rho_B)$ then gives the physical balance equations.
\end{proof}

\section{Boundary ports and thermodynamic balances}
\label{sec:standard_boundary}

\subsection{Boundary power and port variables}
\label{sec:canonical_boundary}

For the diagonal reactor IPHS formulation of
Section~\ref{sec:standard_iphs}, the zero-order reaction
interconnection is pointwise skew-symmetric and therefore generates
no boundary contribution. The boundary power is induced by the
reversible convective interconnection, the mass-diffusion channel,
and the heat-conduction channel.

Using the Green--Stokes identities associated with
$\mathcal J_c$, $\mathcal J_d$, and $\mathcal J_s$, introduce the
endpoint effort vector
\begin{equation}
    \mathfrak e
    :=
    \begin{bmatrix}
        \varphi e_x\\
        \varphi e_s\\
        \varphi R_d e_s\\
        \varphi r_s e_s
    \end{bmatrix}
    =
    \begin{bmatrix}
        -\rho E_c\\
        A_r\\
        v\\
        T\\
        J\\
        -\rho q/T
    \end{bmatrix}.
    \label{eq:modified_effort}
\end{equation}
The first four components correspond to the convective
thermodynamic traces, whereas the last two are associated with
mass diffusion and heat conduction, respectively.

The corresponding symmetric boundary matrix is
\begin{equation}
    P_e(\varphi)
    =
    \begin{bmatrix}
        P_1 & 0 & g_d & 0\\
        0 & 0 & 0 & g_s\\
        g_d^\top & 0 & 0 & 0\\
        0 & g_s & 0 & 0
    \end{bmatrix},
    \qquad
    g_s=\varphi ,
    \label{eq:Pe_block}
\end{equation}
where $g_d=
    \begin{bmatrix}
        0&-1&0
    \end{bmatrix}^T$. On the thermodynamically admissible domain,
\[
    \det P_e(\varphi)
    =
    -\varphi^2
    \neq0,
\]
so that the endpoint metric has full rank.

Direct contraction gives
\begin{equation}
    \frac12
    \mathfrak e^\top
    P_e(\varphi)
    \mathfrak e
    =
    -\rho vE_c
    -
    A_rJ
    -
    q .
    \label{eq:Pe_power}
\end{equation}
Hence the endpoint quadratic form represents the oriented boundary
power density, i.e. the negative of the positive-$z$ physical axial
energy flux.

Because $g_s=\varphi$, the endpoint metric is state dependent.
The following algebraic normalization converts it into a constant
boundary metric without changing the associated power pairing.

\begin{lemma}[Canonical boundary normalization]
\label{lem:boundary_normalization}
Let $N(\varphi)
=
\operatorname{diag}
\left(
1,1,1,\varphi^{-1/2},1,\varphi^{-1/2}
\right)$, $
\overline P_e
=
N(\varphi)^\top P_e(\varphi)N(\varphi)$, $
\widehat{\mathfrak e}
=
N(\varphi)^{-1}\mathfrak e$, then $\overline P_e$ is independent of $\varphi$ and
\[
\mathfrak e^\top P_e(\varphi)\mathfrak e
=
\widehat{\mathfrak e}^{\,\top}
\overline P_e
\widehat{\mathfrak e}.
\]
\end{lemma}

\begin{proof}
The state dependence of $P_e(\varphi)$ occurs only in the
heat-conduction trace pair through $g_s=\varphi$. The two factors
$\varphi^{-1/2}$ reduce this pair to a constant unit coupling, while
the remaining blocks are unchanged. The power identity follows
directly from $\mathfrak e=N(\varphi)\widehat{\mathfrak e}$.
\end{proof}

The normalized endpoint traces admit the following
power-preserving boundary parametrization:
\[
\begin{aligned}
    u_\partial
    &=
    \frac1{\sqrt2}
    \left[
        \widehat{\mathfrak e}(L)
        +
        \widehat{\mathfrak e}(0)
    \right],\\
    y_\partial
    &=
    \frac1{\sqrt2}
    \overline P_e
    \left[
        \widehat{\mathfrak e}(L)
        -
        \widehat{\mathfrak e}(0)
    \right].
\end{aligned}
\]
Since $\overline P_e$ is symmetric,
\begin{equation}
\begin{aligned}
    y_\partial^\top u_\partial
    &=
    \frac12
    \widehat{\mathfrak e}(L)^\top
    \overline P_e
    \widehat{\mathfrak e}(L)
    -
    \frac12
    \widehat{\mathfrak e}(0)^\top
    \overline P_e
    \widehat{\mathfrak e}(0)
    \\
    &=
    -\left[
        \rho vE_c
        +
        A_rJ
        +
        q
    \right]_0^L .
    \label{eq:boundary_supply}
\end{aligned}
\end{equation}
Thus the boundary variables reproduce exactly the physical axial
power exchange of the open reactor.

\subsection{Hamiltonian and entropy balances}
\label{sec:standard_balances}

\begin{theorem}[Boundary energy and entropy balances]
\label{thm:standard_boundary_consistency}

Assume that the reactor state remains in the thermodynamically
admissible domain, with
\[
    M_d\geq0,
    \qquad
    k\geq0,
    \qquad
    \kappa_j\geq0,
\]
and that the reaction regularity condition
\eqref{eq:reaction_dissipation} holds.
Then the diagonal reactor IPHS formulation
\eqref{eq:standard_complete_iphs} satisfies
\begin{equation}
    \dot H
    =
    y_\partial^\top u_\partial
    +
    \int_0^L q_j\,dz
    =
    -
    \left[
        \rho vE_c
        +
        A_rJ
        +
        q
    \right]_0^L
    +
    \int_0^L q_j\,dz .
    \label{eq:standard_energy_balance}
\end{equation}

Its total entropy satisfies
\begin{equation}
\begin{aligned}
    \dot S
    ={}&
    -
    \left[
        \rho sv+\frac qT
    \right]_0^L
    +
    \int_0^L\frac{q_j}{T_j}\,dz
    \\
    &+
    \int_0^L
    \left(
        \sigma_r
        +
        \sigma_d
        +
        \sigma_q
        +
        \sigma_j
    \right)dz ,
\end{aligned}
\label{eq:standard_entropy_balance}
\end{equation}
where
\[
\begin{aligned}
    \sigma_r
    &=
    \frac{A_rr}{T}
    \geq0,
    &
    \sigma_d
    &=
    \frac{M_d}{T}
    (\partial_zA_r)^2
    \geq0,
    \\
    \sigma_q
    &=
    \frac{k}{T^2}
    (\partial_zT)^2
    \geq0,
    &
    \sigma_j
    &=
    \kappa_j
    \frac{(T_j-T)^2}{TT_j}
    \geq0 .
\end{aligned}
\]
\end{theorem}

\begin{proof}
The reaction interconnection $\mathcal J_r$ is pointwise
skew-symmetric and therefore contributes no Hamiltonian power.
Formal skew-adjointness of
$\mathcal J_c$, $\mathcal J_d$, and $\mathcal J_s$
reduces their Hamiltonian contribution to the endpoint pairing
\eqref{eq:Pe_power}. Lemma~\ref{lem:boundary_normalization} and
\eqref{eq:boundary_supply} then give the axial boundary power, while
the distributed jacket port contributes $e_jf_j=q_j$. This yields
\eqref{eq:standard_energy_balance}.

For the entropy balance,
Proposition~\ref{prop:standard_recovery} gives the axial entropy flux
\[
    -
    \left[
        \rho sv+\frac qT
    \right]_0^L
\]
together with the internal contributions
$\sigma_r$, $\sigma_d$, and $\sigma_q$. For the jacket interaction,
\[
    \frac{q_j}{T}
    =
    \frac{q_j}{T_j}
    +
    q_j
    \left(
        \frac1T-\frac1{T_j}
    \right)
    =
    \frac{q_j}{T_j}
    +
    \sigma_j .
\]
Substitution gives
\eqref{eq:standard_entropy_balance}.
\end{proof}

\section{Structure-preserving spatial discretization}
\label{sec:discretization}

This section develops a structure-preserving spatial discretization of
the diagonal $L_c=0$ formulation established in
Sections~\ref{sec:standard_iphs}--\ref{sec:standard_boundary}.

\subsection{Conservative-coordinate PFEM formulation}
\label{sec:pfem_coordinates}

The continuous formulation of Section~\ref{sec:standard_iphs} uses the
Eulerian thermodynamic coordinates
\[
    \xi_s
    :=
    \begin{bmatrix}
        \varphi & \alpha & v & s
    \end{bmatrix}^{\top},
    \qquad
    \varphi=\frac{1}{\rho},
    \qquad
    \alpha=\alpha_A .
\]
For the spatial approximation, introduce the conservative variables
\begin{equation}
    x_c
    :=
    \begin{bmatrix}
        \rho_A & \rho_B & m & \eta
    \end{bmatrix}^{\top}
    =
    \begin{bmatrix}
        \dfrac{\alpha}{\varphi}
        &
        \dfrac{1-\alpha}{\varphi}
        &
        \dfrac{v}{\varphi}
        &
        \dfrac{s}{\varphi}
    \end{bmatrix}^{\top},
    \label{eq:cons_coordinates}
\end{equation}
where
\[
    \rho=\rho_A+\rho_B,
    \qquad
    m:=\rho v,
    \qquad
    \eta:=\rho s .
\]
For $\rho>0$, the transformation $x_c=\Phi(\xi_s)$ is invertible, with
\[
    \varphi=\frac{1}{\rho_A+\rho_B},
    \alpha=\frac{\rho_A}{\rho_A+\rho_B},
    v=\frac{m}{\rho_A+\rho_B},
    s=\frac{\eta}{\rho_A+\rho_B}.
\]

Let $T_\Phi:=\partial\Phi/\partial\xi_s$ denote the Jacobian of the
state transformation. The conjugate efforts satisfy
\[
    \dot x_c=T_\Phi\dot\xi_s,
    \qquad
    e_{\xi_s}=T_\Phi^\top e_c,
\]
and consequently
\begin{equation}
    e_{\xi_s}^{\top}\dot\xi_s
    =
    e_c^{\top}\dot x_c .
    \label{eq:coordinate_power_invariance}
\end{equation}
Thus the conservative-coordinate representation preserves the
Hamiltonian power pairing of the continuous BC-IPHS model.

In conservative coordinates, the Hamiltonian is
\[
    H[x_c]
    =
    \int_0^L
    \mathcal H_c(x_c)\,dz,
    \qquad
    \mathcal H_c
    =
    \frac{m^2}{2\rho}
    +
    \varepsilon(\eta,\rho_A,\rho_B).
\]
Using $  d\varepsilon
    =
    T\,d\eta
    +
    \mu_A\,d\rho_A
    +
    \mu_B\,d\rho_B$, the conjugate efforts are
\[
    e_A^c
    =
    \mu_A-\frac12v^2,
    \qquad
    e_B^c
    =
    \mu_B-\frac12v^2,
    \qquad
    e_m=v,
    \qquad
    e_\eta=T,
\]
so that $e_A^c-e_B^c=A_r$.

We employ a partitioned finite-element approximation
\cite{CardosoRibeiro2018,CardosoRibeiro2021,
KotyczkaMaschkeLefevre2018,Serhani2019,ToledoZucco2024}.
Let $V_h\subset H^1(0,L)$ denote the finite-dimensional approximation
space for the conservative states and their conjugate efforts, and let
$Q_h\subset L^2(0,L)$ contain the associated gradient and flux variables,
with $\partial_zV_h\subset Q_h$. For a basis $\{\phi_i\}_{i=1}^{n}$ of $V_h$, introduce a
positive-weight quadrature rule $\{(z_q,w_q)\}_{q=1}^{N_q}$ and define
\[
    E_{qi}:=\phi_i(z_q),
    \qquad
    G_{qi}:=\partial_z\phi_i(z_q),
    \qquad
    W:=\operatorname{diag}(w_1,\ldots,w_{N_q}).
\]
To distinguish the finite-element mass matrix from the molar mass $M$,
let $M_V:=E^\top WE$, $    K(a)
    :=
    G^\top W\operatorname{diag}(a^q)E$, for any scalar coefficient field $a$. The quadrature rule is assumed to
yield a nonsingular $M_V$. Endpoint terms generated by the partitioned
weak formulation are retained explicitly as boundary contributions.

Let $\mathbf x_h
    :=
    \begin{bmatrix}
        \boldsymbol\rho_A^\top &
        \boldsymbol\rho_B^\top &
        \mathbf m^\top &
        \boldsymbol\eta^\top
    \end{bmatrix}^{\top}$ collect the finite-element coefficients of the conservative state.
The discrete Hamiltonian is evaluated with the same quadrature rule,
\begin{equation}
    H_h(\mathbf x_h)
    :=
    \sum_{q=1}^{N_q}
    w_q\,
    \mathcal H_c
    \left(
        x_{c,h}(z_q)
    \right).
    \label{eq:discrete_hamiltonian}
\end{equation}
Define
\[
    \mathcal M_h
    :=
    \operatorname{diag}(M_V,M_V,M_V,M_V),
    \qquad
    \mathbf e_h
    :=
    \begin{bmatrix}
        (\mathbf e_A^c)^\top &
        (\mathbf e_B^c)^\top &
        \mathbf e_m^\top &
        \mathbf e_\eta^\top
    \end{bmatrix}^{\top}.
\]
The discrete effort vector is the finite-dimensional Riesz
representative of the Hamiltonian gradient,
\begin{equation}
    \mathcal M_h\mathbf e_h
    =
    \nabla_{\mathbf x_h}H_h(\mathbf x_h),
    \label{eq:discrete_riesz_effort}
\end{equation}
and hence
\[
    \dot H_h
    =
    \mathbf e_h^\top
    \mathcal M_h
    \dot{\mathbf x}_h .
\]

All nonlinear thermodynamic and constitutive quantities are evaluated
from the same discrete state, Riesz effort field, and quadrature rule.
In particular,
\[
    \mathbf T^q
    =
    E\mathbf e_\eta,
    \qquad
    \mathbf A_r^q
    =
    E(\mathbf e_A^c-\mathbf e_B^c).
\]
The auxiliary $Q_h$ gradient and flux variables introduced by the
partitioned weak formulation may be eliminated locally after assembly;
the reduced representation below retains the same power pairing and
endpoint traces.

The construction is independent of the polynomial degree of $V_h$.
Polynomial refinement therefore changes only the spatial approximation
accuracy and does not alter the semi-discrete interconnection structure.
A specific Gauss--Lobatto--Legendre (GLL) realization of this
general construction is introduced in Section~\ref{sec:dimensional_benchmark}.

\subsection{Discrete interconnections and boundary ports}
\label{sec:discrete_interconnections}

For the reversible convective channel, the partitioned weak formulation
gives
\begin{equation}
J_{c,h}
=
\begin{bmatrix}
0 & 0 & K(\rho_{A,h}) & 0\\
0 & 0 & K(\rho_{B,h}) & 0\\
-K(\rho_{A,h})^\top
&
-K(\rho_{B,h})^\top
&
K(m_h)-K(m_h)^\top
&
-K(\eta_h)^\top\\
0&0&K(\eta_h)&0
\end{bmatrix}.
\label{eq:discrete_convective_matrix}
\end{equation}
Thus the internal convective interconnection is power preserving, while
the associated external exchange is carried by the endpoint traces
retained in the weak formulation.

For $L_c=0$, the diffusion and thermal constitutive laws are
\[
    J=-M_d\partial_zA_r,
    \qquad
    q=-k\partial_zT .
\]
In conservative coordinates, introduce
\[
    \widehat R_d
    :=
    -\frac JT
    =
    \frac{M_d}{T}\partial_zA_r,
    \qquad
    \widehat R_q
    :=
    -\frac{q}{T^2}
    =
    \frac{k}{T^2}\partial_zT .
\]
At the quadrature points,
\[
    F_d^q
    :=
    \left[
        G(\mathbf e_A^c-\mathbf e_B^c)
    \right]_q,
    \qquad
    F_q^q
    :=
    (G\mathbf e_\eta)_q,
\]
so that
\[
    \widehat R_d^q
    =
    \frac{M_d^q}{T_q}F_d^q,
    \qquad
    \widehat R_q^q
    =
    \frac{k_q}{T_q^2}F_q^q .
\]
Let
\[
    K_d:=K(\widehat R_d),
    \qquad
    K_q:=K(\widehat R_q).
\]
The diffusion and thermal interconnections are
\begin{equation}
    J_{d,h}
    =
    \begin{bmatrix}
        0&0&0&-K_d\\
        0&0&0& K_d\\
        0&0&0&0\\
        K_d^\top&-K_d^\top&0&0
    \end{bmatrix},
    \qquad
    J_{d,h}^\top=-J_{d,h},
    \label{eq:discrete_diffusion_matrix}
\end{equation}
and
\begin{equation}
    J_{q,h}
    =
    \begin{bmatrix}
        0&0&0&0\\
        0&0&0&0\\
        0&0&0&0\\
        0&0&0&K_q^\top-K_q
    \end{bmatrix},
    \qquad
    J_{q,h}^\top=-J_{q,h}.
    \label{eq:discrete_heat_matrix}
\end{equation}
The non-negative coefficients $M_d/T$ and $k/T^2$ generate the
corresponding entropy production.

The reaction and jacket laws are evaluated consistently at the
quadrature points as
\[
    \mathbf r^q
    =
    \mathbf A_r^q
    \odot
    \boldsymbol\ell_r^q,
    \qquad
    \boldsymbol\ell_r^q\geq0,
    \qquad
    \mathbf q_j^q
    =
    \boldsymbol\kappa_j^q
    \odot
    (\mathbf T_j^q-\mathbf T^q).
\]
Their semi-discrete contributions are
\[
    \mathbf g_{r,h}
    =
    \begin{bmatrix}
        -E^\top W\mathbf r^q\\
         E^\top W\mathbf r^q\\
         \mathbf0\\
         E^\top W
         \left(
             \dfrac{
                 \mathbf A_r^q\odot\mathbf r^q
             }{
                 \mathbf T^q
             }
         \right)
    \end{bmatrix},
    \qquad
    \mathbf g_{j,h}
    =
    \begin{bmatrix}
        \mathbf0\\
        \mathbf0\\
        \mathbf0\\
        E^\top W
        \left(
            \dfrac{\mathbf q_j^q}{\mathbf T^q}
        \right)
    \end{bmatrix}.
\]

Collecting the four physical channels gives
\begin{equation}
    \mathcal M_h\dot{\mathbf x}_h
    =
    \left(
        J_{c,h}
        +
        J_{d,h}
        +
        J_{q,h}
    \right)\mathbf e_h
    +
    \mathbf g_{r,h}
    +
    \mathbf g_{j,h}
    +
    \mathbf b_{\partial,h}.
    \label{eq:complete_semidiscrete_system}
\end{equation}

To describe the axial traces, let $f_c$ denote the convective mass-flow
variable, with
\[
    f_c=\rho v=\rho e_m.
\]
The total species fluxes are
\[
    F_A:=\alpha_Af_c+J,
    \qquad
    F_B:=\alpha_Bf_c-J,
\]
while the convective effort is
\[
    E_c
    =
    \alpha_Ae_A^c
    +
    \alpha_Be_B^c
    +
    e_m^2
    +
    se_\eta
    =
    h+\frac12v^2 .
\]

For a general basis of $V_h$, define
\[
    \boldsymbol\phi(z)
    :=
    \begin{bmatrix}
        \phi_1(z)&\cdots&\phi_n(z)
    \end{bmatrix}^{\top},
    \qquad
    B_\partial
    :=
    \begin{bmatrix}
        -\boldsymbol\phi(0)&
        \boldsymbol\phi(L)
    \end{bmatrix}.
\]
For an endpoint-interpolatory nodal basis, including the GLL basis used
in Section~\ref{sec:dimensional_benchmark},
\[
    B_\partial
    =
    \begin{bmatrix}
        -e_1&e_n
    \end{bmatrix}.
\]
For any scalar field $a$, let
\[
    (a)_\partial
    :=
    \begin{bmatrix}
        a(0)\\
        a(L)
    \end{bmatrix}.
\]
The boundary contribution to the semi-discrete state equation is
\begin{equation}
    \mathbf b_{\partial,h}
    =
    \begin{bmatrix}
        -B_\partial(F_A)_\partial\\[1mm]
        -B_\partial(F_B)_\partial\\[1mm]
        -B_\partial(e_mf_c)_\partial\\[1mm]
        -B_\partial
        \left(
            sf_c+\dfrac qT
        \right)_\partial
    \end{bmatrix}.
    \label{eq:discrete_boundary_vector}
\end{equation}
All endpoint quantities are evaluated through the same finite-element
trace operator.

For the physical trace $f_c=\rho v$, direct contraction gives
\begin{equation}
    \mathbf e_h^\top
    \mathbf b_{\partial,h}
    =
    -
    \left[
        \rho vE_c+A_rJ+q
    \right]_0^L .
    \label{eq:discrete_boundary_power_physical}
\end{equation}
Because the canonical normalization introduced in
Section~\ref{sec:standard_boundary} is power preserving, this weak-form
trace is equivalent to the BC-IPHS boundary supply rate,
\[
    \mathcal P_{\partial,h}
    :=
    \mathbf e_h^\top
    \mathbf b_{\partial,h}
    =
    y_{\partial,h}^{\top}u_{\partial,h}.
\]
Likewise, the boundary contribution to the discrete entropy rate is
\[
    \left.\dot S_h\right|_\partial
    =
    -
    \left[
        \rho sv+\frac qT
    \right]_0^L .
\]

For the open-reactor realization, the characteristic count is used only
to identify the missing hydrodynamic boundary datum of the convective
principal subsystem; it is not a well-posedness analysis of the complete
hyperbolic--parabolic reactor problem. The characteristic speeds are
\[
    v-c,\qquad
    v,\qquad
    v,\qquad
    v+c,
    \qquad
    c^2
    =
    \left.
    \frac{\partial p}{\partial\rho}
    \right|_{s,\alpha_A}
    >0 .
\]
Hence, in the positive-subsonic regime $0<v<c$, three convective
characteristics enter at $z=0$ and one enters at $z=L$.

The species and thermal boundary data are imposed separately through
total flux conditions. The inlet conditions are
\begin{align}
    (\rho_Av+J)(0,t)
    &=
    F_A^{\mathrm{in}}(t),
    \nonumber\\
    (\rho_Bv-J)(0,t)
    &=
    F_B^{\mathrm{in}}(t),
    \nonumber\\
    (\rho vE_c+A_rJ+q)(0,t)
    &=
    F_H^{\mathrm{in}}(t),
    \label{eq:physical_inlet_conditions}
\end{align}
whereas the outlet conditions are
\begin{equation}
    J(L,t)=0,
    \qquad
    q(L,t)=0,
    \qquad
    p(L,t)=p_{\mathrm{out}}(t).
    \label{eq:physical_outlet_conditions}
\end{equation}
Since $F_A+F_B=\rho v$, the two inlet species fluxes determine the inlet
convective mass flow; no additional inlet-velocity condition is imposed.
The outlet conditions $J=0$ and $q=0$ are retained in physical flux form
and are not replaced by independent zero-gradient conditions.

In the numerical realization, the outlet pressure is enforced through
the existing convective boundary port. Write
\[
    \mathcal M_h\dot{\mathbf x}_h
    =
    \mathbf F_h(\mathbf x_h,t)
    +
    \mathbf b_L(\mathbf x_h)f_{c,L},
\]
together with the algebraic constraint
\begin{equation}
    g(\mathbf x_h,t)
    :=
    p_h(L,t)-p_{\mathrm{out}}(t)
    =
    0.
    \label{eq:boundary_DAE_constraint}
\end{equation}
Differentiating the constraint gives
\[
    g_{\mathbf x_h}
    \mathcal M_h^{-1}
    \left(
        \mathbf F_h+\mathbf b_Lf_{c,L}
    \right)
    +
    g_t
    =
    0 .
\]
Whenever
\[
    g_{\mathbf x_h}
    \mathcal M_h^{-1}
    \mathbf b_L
    \neq0,
\]
the outlet convective flow is
\begin{equation}
    f_{c,L}
    =
    -
    \frac{
        g_{\mathbf x_h}
        \mathcal M_h^{-1}\mathbf F_h
        +
        g_t
    }{
        g_{\mathbf x_h}
        \mathcal M_h^{-1}\mathbf b_L
    }.
    \label{eq:index_reduced_outlet_flow}
\end{equation}
This index reduction only enforces the prescribed outlet-pressure
constraint. It introduces neither an additional dissipative mechanism
nor a new boundary power channel.

\subsection{Discrete energy and entropy balances}
\label{sec:discrete_thermodynamics}

Let $\mathbf c$ denote the coefficient vector representing the constant
function $1$ in $V_h$, so that
\[
    E\mathbf c=\mathbf1,
    \qquad
    G\mathbf c=\mathbf0.
\]
The discrete total entropy is
\begin{equation}
    S_h
    :=
    \mathbf c^\top M_V\boldsymbol\eta .
    \label{eq:discrete_entropy}
\end{equation}

\begin{theorem}[Discrete energy and entropy balances]
\label{thm:discrete_thermodynamic_consistency}

Consider the semi-discrete system
\eqref{eq:complete_semidiscrete_system}.
Assume that the quadrature weights are positive, $M_V$ is nonsingular,
$1\in V_h$, and, at every quadrature point,
\[
    \rho_q>0,\qquad
    T_q>0,\qquad
    T_{j,q}>0,\qquad
    M_{d,q},\,k_q,\,\ell_{r,q},\,\kappa_{j,q}\geq0.
\]
Assume also that all discrete thermodynamic and constitutive quantities
are evaluated using the same quadrature rule and Riesz effort field,
with endpoint quantities evaluated through the same finite-element
trace operator.

Then the spatially semi-discrete system satisfies
\begin{equation}
    \dot H_h
    =
    \mathcal P_{\partial,h}+P_{j,h},
    \qquad
    P_{j,h}:=\sum_q w_q q_j^q,
    \label{eq:discrete_first_law}
\end{equation}
where
\[
    \mathcal P_{\partial,h}
    =
    -\left[E_cf_c+A_rJ+q\right]_0^L.
\]
For the physical trace $f_c=\rho v$,
\begin{equation}
    \mathcal P_{\partial,h}
    =
    y_{\partial,h}^\top u_{\partial,h}
    =
    -\left[\rho vE_c+A_rJ+q\right]_0^L .
    \label{eq:discrete_physical_power}
\end{equation}

Furthermore,
\begin{equation}
\begin{aligned}
    \dot S_h
    ={}&
    \mathcal F_{\partial,h}^{S}
    +\Phi_{j,h}
    +\Sigma_{r,h}
    +\Sigma_{d,h}
    +\Sigma_{q,h}
    +\Sigma_{j,h},
\end{aligned}
\label{eq:discrete_second_law}
\end{equation}
where
\[
    \mathcal F_{\partial,h}^{S}
    :=
    -\left[sf_c+\frac qT\right]_0^L,
    \qquad
    \Phi_{j,h}
    :=
    \sum_q w_q\frac{q_j^q}{T_j^q},
\]
and
\begin{equation}
\begin{aligned}
    \Sigma_{r,h}
    &:=
    \sum_q w_q
    \frac{(A_{r,q})^2\ell_{r,q}}{T_q}
    \geq0,
    \\
    \Sigma_{d,h}
    &:=
    \sum_q w_q
    \frac{M_{d,q}}{T_q}(F_d^q)^2
    \geq0,
    \\
    \Sigma_{q,h}
    &:=
    \sum_q w_q
    \frac{k_q}{T_q^2}(F_q^q)^2
    \geq0,
    \\
    \Sigma_{j,h}
    &:=
    \sum_q w_q\kappa_{j,q}
    \frac{(T_j^q-T_q)^2}{T_qT_j^q}
    \geq0.
\end{aligned}
\label{eq:discrete_entropy_production}
\end{equation}
For $f_c=\rho v$,
\[
    \mathcal F_{\partial,h}^{S}
    =
    -\left[\rho sv+\frac qT\right]_0^L .
\]

Hence the spatial semi-discretization satisfies the discrete counterpart
of the continuous Hamiltonian balance and an entropy balance with
non-negative reaction, diffusion, heat-conduction, and jacket
entropy-production contributions.
\end{theorem}

\begin{proof}
From the Riesz representation,
\[
    \dot H_h
    =
    \mathbf e_h^\top
    \mathcal M_h\dot{\mathbf x}_h .
\]
Substitution of \eqref{eq:complete_semidiscrete_system}, together with
\[
    J_{c,h}^\top=-J_{c,h},
    \qquad
    J_{d,h}^\top=-J_{d,h},
    \qquad
    J_{q,h}^\top=-J_{q,h},
\]
eliminates all internal interconnection powers. The reaction terms cancel
pointwise since $A_r=e_A^c-e_B^c$:
\[
    -e_A^{c,q}r_q
    +e_B^{c,q}r_q
    +T_q\frac{A_{r,q}r_q}{T_q}
    =0.
\]
The jacket contribution is
\[
    \mathbf e_h^\top\mathbf g_{j,h}
    =
    \sum_q w_q q_j^q
    =
    P_{j,h}.
\]
Together with the boundary contribution, this proves \eqref{eq:discrete_first_law}.

For the entropy balance, premultiplying the $\eta$ equation by
$\mathbf c^\top$ and using
$E\mathbf c=\mathbf1$ and $G\mathbf c=\mathbf0$
reduces the discrete divergence terms to
$\mathcal F_{\partial,h}^{S}$.
The reaction, diffusion, and heat-conduction terms yield, respectively,
$\Sigma_{r,h}$, $\Sigma_{d,h}$, and $\Sigma_{q,h}$ as defined in
\eqref{eq:discrete_entropy_production}.
Finally,
\[
    \frac{q_j^q}{T_q}
    =
    \frac{q_j^q}{T_j^q}
    +
    \kappa_{j,q}
    \frac{(T_j^q-T_q)^2}{T_qT_j^q},
    \qquad
    q_j^q=\kappa_{j,q}(T_j^q-T_q),
\]
so that summation with the same positive quadrature rule gives
$\Phi_{j,h}+\Sigma_{j,h}$ and proves
\eqref{eq:discrete_second_law}.
\end{proof}

\begin{remark}
Theorem~\ref{thm:discrete_thermodynamic_consistency} concerns the
spatially semi-discrete model, with time left continuous; no exact
Hamiltonian or entropy identity is claimed for the subsequent time
integration. The outlet-pressure index reduction only determines the
existing convective boundary flow and introduces no additional power
channel or entropy-production mechanism.
\end{remark}

\section{Numerical results}
\subsection{Standard tubular-reactor model}
\label{sec:dimensional_benchmark}

The numerical study considers a dimensional open tubular reactor under the
jacket-temperature operating scenario adopted in our previous reactor
study~\cite{zhou2021thermodynamic}. The numerical case corresponds to the
diagonal transport model developed in
Sections~\ref{sec:standard_iphs}--\ref{sec:discretization}, with
$L_c=0$. The non-diagonal Soret--Dufour extension is considered separately
in Section~\ref{sec:numerical_thermodynamics}.

The thermodynamic and kinetic quantities are expressed on the mass-specific
basis introduced in Section~\ref{sec:thermodynamics}. Equal molar masses are
used,
\[
    M_A=M_B=M=0.5~\mathrm{kg\,mol^{-1}},
    \qquad
    R_s=\frac{R}{M}
    =16.628~\mathrm{J\,kg^{-1}K^{-1}}.
\]
For $i\in\{A,B\}$, the caloric functions are
\[
    h_i(T)
    =
    h_{i,\mathrm{ref}}
    +
    c_{p,i}(T-T_{\mathrm{ref}}),
    \qquad
    s_i^0(T)
    =
    s_{i,\mathrm{ref}}
    +
    c_{p,i}\ln\!\left(\frac{T}{T_{\mathrm{ref}}}\right),
\]
and the specific mixture entropy is
\[
    s
    =
    \alpha_A s_A^0(T)
    +
    \alpha_B s_B^0(T)
    -
    R_s
    \left(
        \alpha_A\ln\alpha_A
        +
        \alpha_B\ln\alpha_B
    \right).
\]

A weakly compressible liquid closure is employed to retain the distributed
momentum dynamics while keeping the reactor in a very-low-Mach regime,
\begin{equation}
    p
    =
    p_{\mathrm{ref}}
    +
    K_{\mathrm{bulk}}
    \left(
        \frac{\rho}{\rho_{\mathrm{ref}}}-1
    \right).
    \label{eq:num_pressure_closure}
\end{equation}
To generate this pressure law from the same thermodynamic potential used in
the Hamiltonian, the specific compression Helmholtz contribution is chosen as
\[
    a_c(\rho)
    =
    \frac{K_{\mathrm{bulk}}}{\rho_{\mathrm{ref}}}
    \ln\!\left(\frac{\rho}{\rho_{\mathrm{ref}}}\right)
    +
    (p_{\mathrm{ref}}-K_{\mathrm{bulk}})
    \left(
        \frac{1}{\rho_{\mathrm{ref}}}
        -
        \frac{1}{\rho}
    \right),
    \qquad
    \rho^2\frac{da_c}{d\rho}=p.
\]
Defining
\[
    a_i^0(T)
    :=
    h_i(T)
    -
    \frac{p_{\mathrm{ref}}}{\rho_{\mathrm{ref}}}
    -
    Ts_i^0(T),
\]
the chemical potentials are evaluated as
\begin{equation}
    \mu_i
    =
    a_i^0(T)
    +
    R_sT\ln\alpha_i
    +
    a_c(\rho)
    +
    \frac{p}{\rho}.
    \label{eq:num_chemical_potential}
\end{equation}
The common compression contribution cancels from
$A_r=\mu_A-\mu_B$, consistently with the thermodynamic specialization of
Section~\ref{sec:thermodynamics}.

The reaction law and its regular factorization are
\begin{equation}
\begin{aligned}
    r
    &=
    k_f(T)\rho_A
    \left[
        1-
        \exp\!\left(-\frac{A_r}{R_sT}\right)
    \right],
    k_f(T)
    =
    k_0\exp\!\left(-\frac{E_a}{RT}\right),
    \\
    r
    &=
    A_r\ell_r(A_r),
    \ell_r(A_r)
    =
    k_f(T)\rho_A
    \frac{
        1-\exp[-A_r/(R_sT)]
    }{
        A_r
    }
    \geq0,
     A_r\neq0,
\end{aligned}
\label{eq:numerical_reaction_law}
\end{equation}
with the continuous limit
\[
    \ell_r(0)
    =
    \frac{k_f(T)\rho_A}{R_sT}
    >0.
\]
Hence the numerical reaction law satisfies the reaction-channel regularity
condition of Section~\ref{sec:standard_interconnections}.

For the diagonal model, the constitutive transport laws are
\begin{equation}
    M_d
    =
    \frac{D\rho\alpha_A\alpha_B}{R_sT},
    \qquad
    J=-M_d\partial_zA_r,
    \qquad
    q=-k\partial_zT .
    \label{eq:num_standard_transport}
\end{equation}

Table~\ref{tab:benchmark_parameters} summarizes the dimensional parameters.
Parameters marked as present-benchmark quantities are numerical
specializations and are not claimed to have been experimentally identified
in~\cite{zhou2021thermodynamic}.

\begin{table}[t]
    \centering
    \caption{Dimensional parameters used in the numerical benchmark.}
    \label{tab:benchmark_parameters}

    \footnotesize
    \setlength{\tabcolsep}{3pt}
    \renewcommand{\arraystretch}{1.10}

    \begin{tabularx}{\columnwidth}{
        @{}
        >{\raggedright\arraybackslash}p{0.25\columnwidth}
        >{\raggedright\arraybackslash}p{0.22\columnwidth}
        >{\raggedright\arraybackslash}X
        @{}
    }
        \hline
        Quantity & Symbol & Value \\
        \hline

        Reactor geometry
        &
        $L,\;V,\;d$
        &
        $1~\mathrm{m},\;
        10^{-3}~\mathrm{m^3},\;
        3.568\times10^{-2}~\mathrm{m}$ ${}^{a,b,c}$
        \\

        Molar mass
        &
        $M_A=M_B$
        &
        $0.5~\mathrm{kg\,mol^{-1}}$ ${}^{c}$
        \\

        Reference state
        &
        $T_{\mathrm{ref}},\rho_{\mathrm{ref}},p_{\mathrm{ref}}$
        &
        $300~\mathrm{K},\;
        1000~\mathrm{kg\,m^{-3}},\;
        101325~\mathrm{Pa}$ ${}^{a,b,c}$
        \\

        Heat capacities
        &
        $c_{p,A},c_{p,B}$
        &
        $300.96,\;240.00~
        \mathrm{J\,kg^{-1}K^{-1}}$ ${}^{a}$
        \\

        Reference enthalpies
        &
        $h_{A,\mathrm{ref}},h_{B,\mathrm{ref}}$
        &
        $0,\;-18.30~
        \mathrm{kJ\,kg^{-1}}$ ${}^{a}$
        \\

        Reference entropies
        &
        $s_{A,\mathrm{ref}},s_{B,\mathrm{ref}}$
        &
        $420.8,\;360.4~
        \mathrm{J\,kg^{-1}K^{-1}}$ ${}^{a}$
        \\

        Compressibility
        &
        $K_{\mathrm{bulk}}$
        &
        $2.2~\mathrm{GPa}$ ${}^{c}$
        \\

        Reaction kinetics
        &
        $k_0,\;E_a$
        &
        $1.2\times10^9~\mathrm{s^{-1}},\;
        72.335~\mathrm{kJ\,mol^{-1}}$ ${}^{a}$
        \\

        Axial transport
        &
        $D,\;k$
        &
        $4.5\times10^{-5}~\mathrm{m^2\,s^{-1}},\;
        12.1716~\mathrm{W\,m^{-1}K^{-1}}$ ${}^{a,c}$
        \\

        Jacket exchange
        &
        $U_j,\;\kappa_j$
        &
        $100~\mathrm{W\,m^{-2}K^{-1}},\;
        1.1210\times10^4~
        \mathrm{W\,m^{-3}K^{-1}}$ ${}^{b,c}$
        \\

        Operating condition
        &
        $T_{\mathrm{in}},T_j,$
        $\alpha_{A,\mathrm{in}},v_{\mathrm{in}}$
        &
        $330~\mathrm{K},\;
        350\!\rightarrow\!370~\mathrm{K},\;
        0.9999,\;
        5\times10^{-3}~\mathrm{m\,s^{-1}}$
        ${}^{a,c}$
        \\

        \hline
    \end{tabularx}

    \vspace{1mm}

    \begin{minipage}{\columnwidth}
        \scriptsize
        ${}^{a}$ Inherited from, or converted to the mass-specific SI
        basis from, the reactor data in
        \cite{zhou2021thermodynamic}.
        ${}^{b}$ Derived from the listed geometry or reference data.
        ${}^{c}$ Specified for the present benchmark.
    \end{minipage}
\end{table}

The nominal inlet convective mass flux is
\[
    f_{c,\mathrm{in}}
    =
    \rho_{\mathrm{ref}}v_{\mathrm{in}}.
\]
The prescribed species and energy fluxes are
\[
    F_A^{\mathrm{in}}
    =
    \alpha_{A,\mathrm{in}}f_{c,\mathrm{in}},
    \qquad
    F_B^{\mathrm{in}}
    =
    \alpha_{B,\mathrm{in}}f_{c,\mathrm{in}},
    \qquad
    F_H^{\mathrm{in}}
    =
    E_{c,\mathrm{in}}f_{c,\mathrm{in}},
\]
where
$\alpha_{B,\mathrm{in}}=1-\alpha_{A,\mathrm{in}}$.
These are total-flux boundary conditions; no additional independent inlet
conditions are imposed on $J$ or $q$. At the outlet,
\[
    J(L,t)=0,
    \qquad
    q(L,t)=0,
    \qquad
    p(L,t)=p_{\mathrm{out}},
\]
with the pressure condition enforced by the index-reduced boundary DAE of
Section~\ref{sec:discrete_interconnections}.

The reactor is initialized at the nonuniform steady state corresponding to
$T_j=350~\mathrm{K}$. At $t=0$, the jacket temperature is changed from
$350$ to $370~\mathrm{K}$, while the inlet fluxes and outlet pressure remain
unchanged. Since the jacket exchanges heat along the reactor length, this
perturbation acts as a distributed thermal forcing rather than as an
inlet-temperature disturbance.

The spatial approximation uses a single degree-$60$ nodal Lagrange
element with Gauss--Lobatto--Legendre (GLL) nodes, yielding $61$ global
nodes. This provides a specific polynomial implementation of the PFEM
construction in Section~\ref{sec:dimensional_benchmark}; no artificial viscosity,
filtering, stabilization, penalty term, pressure correction, or empirical
numerical flux is introduced.

As an independent spatial consistency check, the diagonal ($L_c=0$)
semi-discrete initial state was compared with the corresponding continuous
steady boundary-value solution, giving
\[
    \max_z|T_h-T_{\mathrm{BVP}}|
    \simeq4.0\times10^{-7}~\mathrm{K},
    \qquad
    \max_z|\alpha_{B,h}-\alpha_{B,\mathrm{BVP}}|
    \simeq5.8\times10^{-8}.
\]
These discrepancies provide an independent check of the spatial
realization and are not part of the structure-preserving proof.

\begin{figure}[t]
    \centering
    \includegraphics[width=\columnwidth]
    {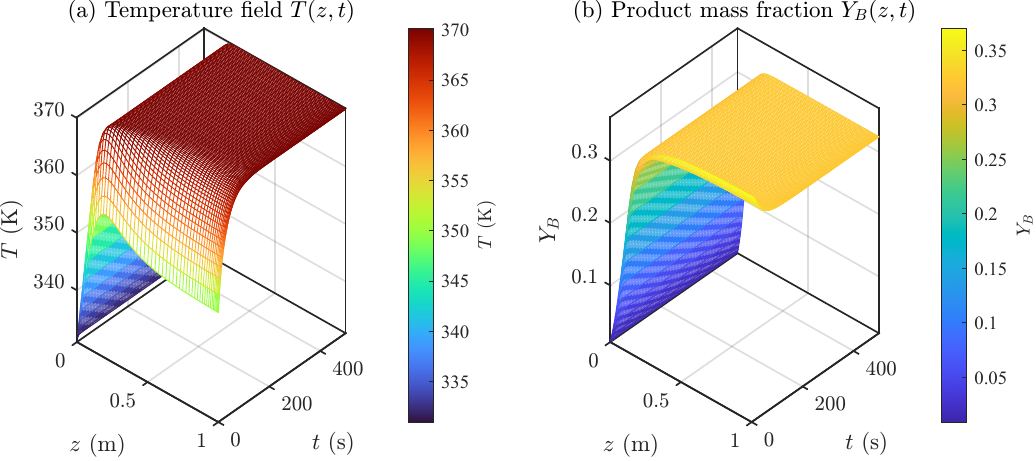}
    \caption{Distributed response of the diagonal benchmark following the
    jacket-temperature step:
    (a) temperature field $T(z,t)$ and
    (b) product mass-fraction field $\alpha_B(z,t)$.
    The initial condition is the steady state corresponding to
    $T_j=350~\mathrm{K}$.}
    \label{fig:space_time_dynamics}
\end{figure}

Figure~\ref{fig:space_time_dynamics} shows the distributed temperature and
composition response over $0\leq z\leq L$ and
$0\leq t\leq500~\mathrm{s}$.

The jacket-temperature increase raises the distributed heat input
$q_j=\kappa_j(T_j-T)$, while the inlet continues to supply fluid at
$T_{\mathrm{in}}=330~\mathrm{K}$. The transient therefore reflects the
competition among distributed wall heating, axial convection, and heat
conduction rather than an inlet-generated thermal wave. The composition
response is governed jointly by the temperature dependence of $k_f(T)$
and $A_r/(R_sT)$, convective residence, and affinity-driven diffusion.

\begin{figure}[t]
    \centering
    \includegraphics[width=\columnwidth]
    {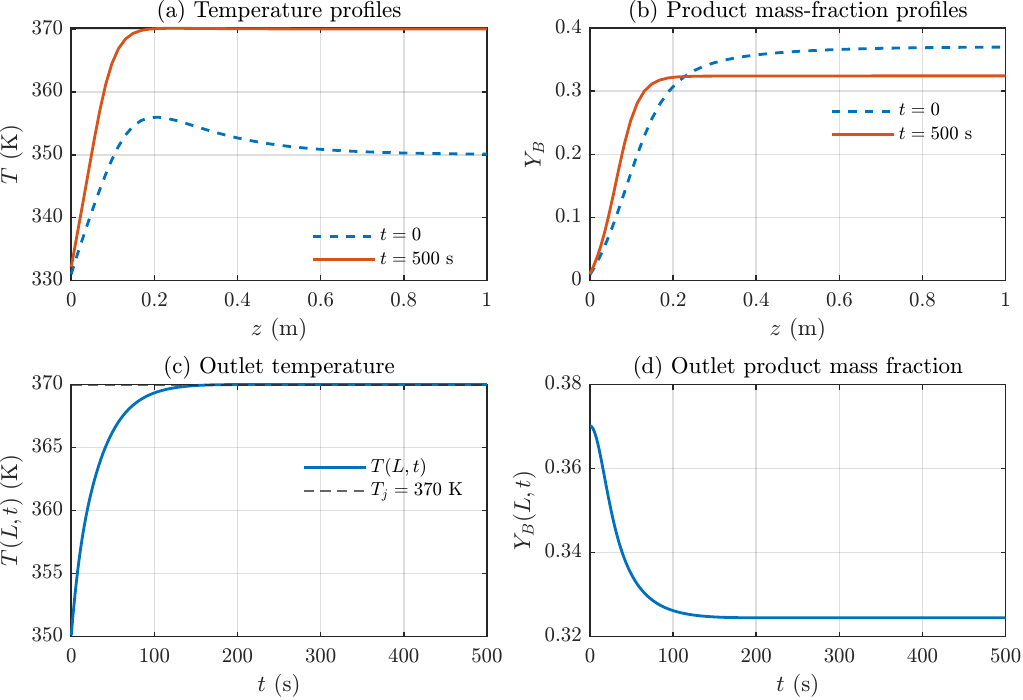}
    \caption{Diagonal-model spatial profiles and outlet response:
    (a) initial and final temperature profiles;
    (b) initial and final product mass-fraction profiles;
    (c) outlet temperature $T(L,t)$; and
    (d) outlet product mass fraction $\alpha_B(L,t)$.}
    \label{fig:profiles_outlet}
\end{figure}

The initial and final spatial profiles and the corresponding outlet traces
are shown in Fig.~\ref{fig:profiles_outlet}.

Figure~\ref{fig:profiles_outlet} makes the steady-profile transition
explicit. The new temperature profile contains an upstream development
region maintained by the continuous $330~\mathrm{K}$ feed and a nearly
isothermal downstream region close to $370~\mathrm{K}$. The product
fraction develops from its near-zero inlet value and approaches a lower
downstream level. Thus the jacket step transfers the reactor between two
distributed open-system steady states rather than merely changing the
outlet values. At $t=500~\mathrm{s}$,
\[
    T(L,500)\simeq369.9997~\mathrm{K},
    \qquad
    \alpha_B(L,500)\simeq0.324487.
\]

The momentum dynamics remain in the intended weakly compressible,
very-low-Mach regime. From \eqref{eq:num_pressure_closure},
\[
    c
    =
    \sqrt{\frac{K_{\mathrm{bulk}}}{\rho_{\mathrm{ref}}}}
    \simeq1.48\times10^3~\mathrm{m\,s^{-1}},
    \qquad
    Ma_{\mathrm{in}}
    \simeq3.37\times10^{-6}.
\]
Over the jacket-step transient, the pressure deviation remains below
$3\times10^{-1}~\mathrm{Pa}$ and the velocity deviation below
$2\times10^{-8}~\mathrm{m\,s^{-1}}$. The observable response is therefore
dominated by thermal, reaction, and species-transport dynamics, while the
distributed momentum state remains part of the Eulerian convection and
boundary-energy structure.

The weak velocity response to the jacket-temperature step does not imply
that the convection velocity is prescribed. The transported momentum
density $m=\rho v$ is a dynamic state, with
\[
    v(z,t)
    =
    \frac{m(z,t)}
    {\rho_A(z,t)+\rho_B(z,t)}.
\]
To expose this dynamic degree of freedom independently of thermal forcing,
an additional feed-throughput transition is applied with
$T_j=350~\mathrm{K}$ and unchanged outlet pressure. The inlet convective
mass flux is increased by $5\%$ through a $10~\mathrm{s}$ linear ramp,
\begin{equation}
    f_{c,\mathrm{in}}(t)
    =
    \beta(t)f_{c,0},
    \beta(t)
    =
    \begin{cases}
        1+0.05\,t/10, & 0\leq t<10~\mathrm{s},\\
        1.05, & t\geq10~\mathrm{s}.
    \end{cases}
    \label{eq:throughput_ramp}
\end{equation}
The feed composition and specific inlet-energy level are unchanged, so the
prescribed extensive inlet fluxes are scaled consistently:
\begin{equation}
    F_A^{\mathrm{in}}(t)
    =
    \beta(t)F_{A,0}^{\mathrm{in}},
    F_B^{\mathrm{in}}(t)
    =
    \beta(t)F_{B,0}^{\mathrm{in}},
    F_H^{\mathrm{in}}(t)
    =
    \beta(t)F_{H,0}^{\mathrm{in}}.
    \label{eq:throughput_scaled_fluxes}
\end{equation}
Thus no independent inlet-velocity boundary condition is introduced; the
velocity response follows from the momentum dynamics and the existing
flux/pressure boundary ports.

\begin{figure}[t]
    \centering
    \includegraphics[width=\columnwidth]
    {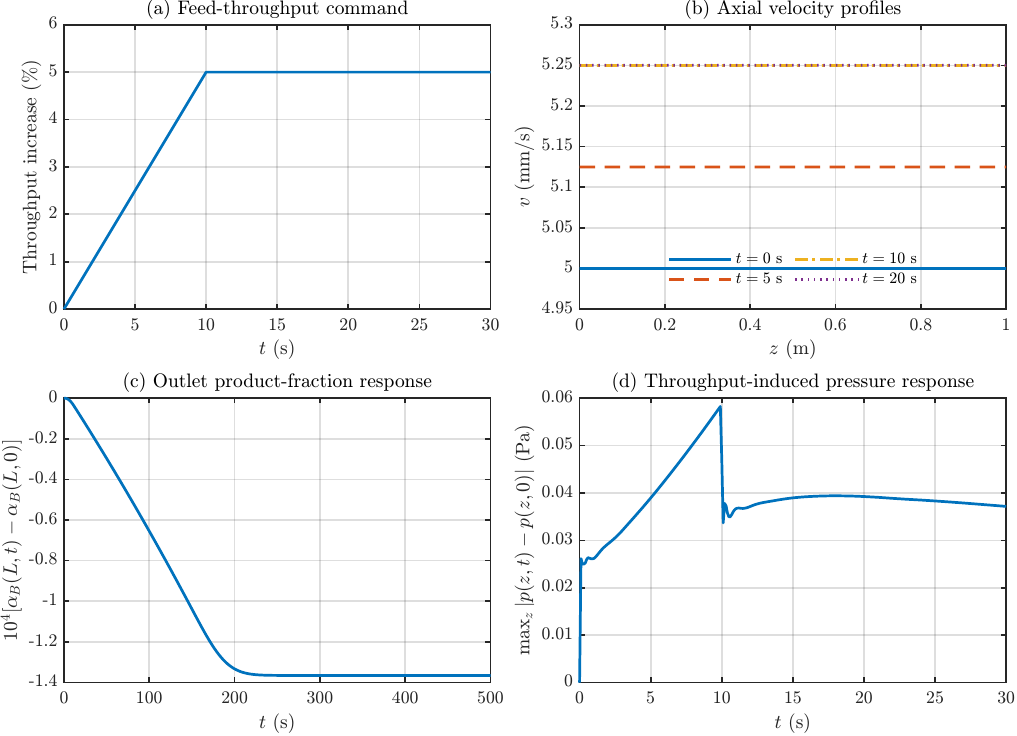}
    \caption{Dynamic convective response to a controlled $5\%$
    feed-throughput increase applied over $10~\mathrm{s}$:
    (a) prescribed inlet convective mass flux;
    (b) distributed velocity profiles at selected times;
    (c) change in outlet product mass fraction relative to the initial
    steady state; and
    (d) maximum pressure-field deviation from the initial profile.
    The jacket temperature is held at $350~\mathrm{K}$ and the outlet
    pressure remains prescribed throughout the transition.}
    \label{fig:throughput_velocity}
\end{figure}

Figure~\ref{fig:throughput_velocity} shows that the mean convection
velocity increases from approximately
$5.00\times10^{-3}$ to
$5.25\times10^{-3}~\mathrm{m\,s^{-1}}$, while remaining nearly uniform in
space; its maximum instantaneous spatial spread is only
$2.66\times10^{-8}~\mathrm{m\,s^{-1}}$. This is consistent with the high
bulk modulus and very-low-Mach regime, for which the hydrodynamic
adjustment is much faster than the $10~\mathrm{s}$ feed-flow ramp.
Accordingly, the convective residence-time scale $L/\bar v$ decreases
from approximately $200$ to $190.5~\mathrm{s}$.

The shorter residence time produces a small but systematic reduction in
the outlet product fraction,
\[
    \alpha_B(L,500)-\alpha_B(L,0)
    \simeq-1.37\times10^{-4},
\]
while
\[
    \max_{z,t}|p(z,t)-p(z,0)|
    \simeq5.82\times10^{-2}~\mathrm{Pa}.
\]
The latter is not interpreted as an industrial pressure-drop prediction,
since distributed wall friction is not included in the present momentum
model. Rather, this benchmark demonstrates that the convection velocity
is a dynamic distributed variable and that a throughput change propagates
through the momentum state to the reactor residence and composition
dynamics.

\subsection{Energy--entropy balances and Soret--Dufour transport}
\label{sec:numerical_thermodynamics}

\begin{figure}[t]
    \centering
    \includegraphics[width=\columnwidth]
    {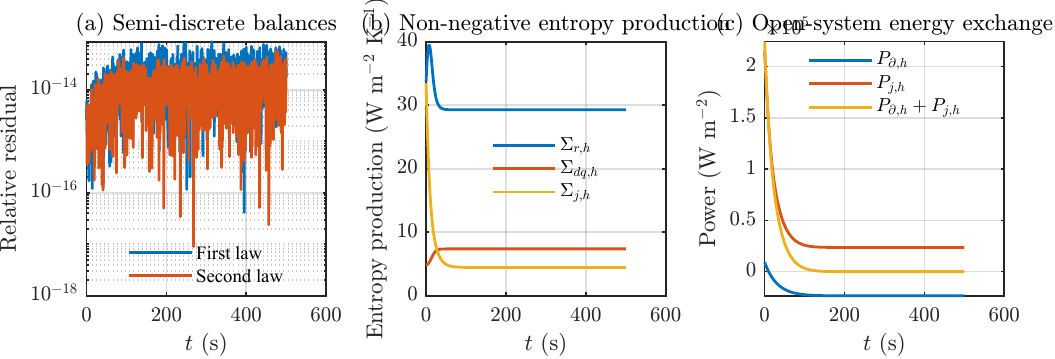}
    \caption{Hamiltonian and entropy balances of the diagonal semi-discrete model:
    (a) normalized residuals of the Hamiltonian and entropy balances;
    (b) non-negative entropy-production contributions due to reaction,
    diffusion, heat conduction, and jacket heat exchange; and
    (c) axial-boundary and jacket contributions to the Hamiltonian balance.}
    \label{fig:thermo_structure}
\end{figure}

Figure~\ref{fig:thermo_structure} evaluates the semi-discrete Hamiltonian
and entropy identities of
Theorem~\ref{thm:discrete_thermodynamic_consistency}
along the diagonal $L_c=0$ trajectory. For clarity, define the normalized
balance residuals
\begin{equation}
\begin{aligned}
\varepsilon_H
&=
\frac{
    \left|
        \dot H_h-\mathcal P_{\partial,h}-P_{j,h}
    \right|
}{
    \max\!\left\{
        1,\,
        |\dot H_h|
        +|\mathcal P_{\partial,h}|
        +|P_{j,h}|
    \right\}
},\\[1mm]
\varepsilon_S
&=
\frac{
    \left|
        \dot S_h
        -\mathcal F_{\partial,h}^{S}
        -\Phi_{j,h}
        -\Sigma_{r,h}
        -\Sigma_{d,h}
        -\Sigma_{q,h}
        -\Sigma_{j,h}
    \right|
}{
    \max\!\left\{
        1,\,
        |\dot S_h|
        +|\mathcal F_{\partial,h}^{S}|
        +|\Phi_{j,h}|
        +|\Sigma_{r,h}|
        +|\Sigma_{d,h}|
        +|\Sigma_{q,h}|
        +|\Sigma_{j,h}|
    \right\}
}.
\end{aligned}
\label{eq:numerical_balance_residuals}
\end{equation}

The first- and second-law residuals remain at numerical round-off level,
while
\[
    \Sigma_{r,h}\geq0,
    \qquad
    \Sigma_{d,h}\geq0,
    \qquad
    \Sigma_{q,h}\geq0,
    \qquad
    \Sigma_{j,h}\geq0
\]
throughout the transient. These quantities represent entropy production
due to chemical reaction, affinity-driven mass diffusion, heat conduction,
and finite-temperature-difference jacket exchange, respectively.

The power and entropy diagnostics recover the identities of
Theorem~\ref{thm:discrete_thermodynamic_consistency}. In particular, as
the reactor approaches the new steady profile, $\dot H_h\rightarrow0$
does not imply that the jacket power vanishes; rather,
\[
    \mathcal P_{\partial,h}+P_{j,h}\rightarrow0,
\]
so that the distributed jacket input is balanced by the net axial
boundary-energy flux, as expected for an open reactor at steady operation.

The numerical trajectory therefore corroborates the spatial semi-discrete
thermodynamic identities. These diagnostics refer to the time-continuous semi-discrete PFEM system.
The subsequent \texttt{ode15s} integration is not a
structure-preserving time discretization, and no exact fully discrete
Hamiltonian or entropy identity is claimed. Structure-preserving temporal
discretization constitutes a separate issue and is not considered here;
see, e.g., \cite{KotyczkaLefevre2019}.

The preceding structural and numerical results concern the diagonal
transport model. To quantify the effect of reciprocal diffusion--heat
coupling, an additional constitutive model is considered with
\[
    L_c=\chi\sqrt{M_dkT^3},
    \qquad
    \chi=0.10.
\]
Since $L_c^2=\chi^2M_dkT^3$, the corresponding Onsager matrix remains
positive semidefinite for $|\chi|\leq1$.

At each quadrature point $\ell$, define
\[
    F_{d,\ell}
    =
    \left[G(e_A^c-e_B^c)\right]_\ell,
    \qquad
    F_{q,\ell}
    =
    (Ge_\eta)_\ell,
\]
with thermodynamic forces
\[
    X_{d,\ell}
    =
    -\frac{F_{d,\ell}}{T_\ell},
    \qquad
    X_{q,\ell}
    =
    -\frac{F_{q,\ell}}{T_\ell^2}.
\]
The coupled constitutive law is
\begin{equation}
    \begin{bmatrix}
        J_\ell\\[1mm]
        q_\ell
    \end{bmatrix}
    =
    L_\ell
    \begin{bmatrix}
        X_{d,\ell}\\[1mm]
        X_{q,\ell}
    \end{bmatrix},
    \qquad
    L_\ell
    =
    \begin{bmatrix}
        M_{d,\ell}T_\ell & L_{c,\ell}\\
        L_{c,\ell} & k_\ell T_\ell^2
    \end{bmatrix}.
    \label{eq:numerical_nondiagonal_constitutive}
\end{equation}
The corresponding discrete diffusion--heat entropy production is
\begin{equation}
    \Sigma_{dq,h}
    =
    \sum_{\ell}w_\ell
    \begin{bmatrix}
        X_{d,\ell} & X_{q,\ell}
    \end{bmatrix}
    L_\ell
    \begin{bmatrix}
        X_{d,\ell}\\
        X_{q,\ell}
    \end{bmatrix}
    \geq0.
    \label{eq:numerical_coupled_entropy}
\end{equation}
This non-diagonal constitutive extension is used only for numerical
comparison and is not part of the diagonal structural result established
in Sections~\ref{sec:standard_iphs}--\ref{sec:discretization}.

The two constitutive models are initialized from their respective
$T_j=350~\mathrm{K}$ steady states before the same
$350\rightarrow370~\mathrm{K}$ jacket step is applied. To separate the
transient coupling effect from differences between the initial operating
profiles, define, for $\psi\in\{T,\alpha_B\}$,
\begin{equation}
    \Delta\psi_{\mathrm{step}}(z,t)
    :=
    \left[
        \psi_{0.10}(z,t)-\psi_{0.10}(z,0)
    \right]
    -
    \left[
        \psi_{0}(z,t)-\psi_{0}(z,0)
    \right].
    \label{eq:coupling_step_difference}
\end{equation}

\begin{figure}[t]
    \centering
    \includegraphics[width=\columnwidth]
    {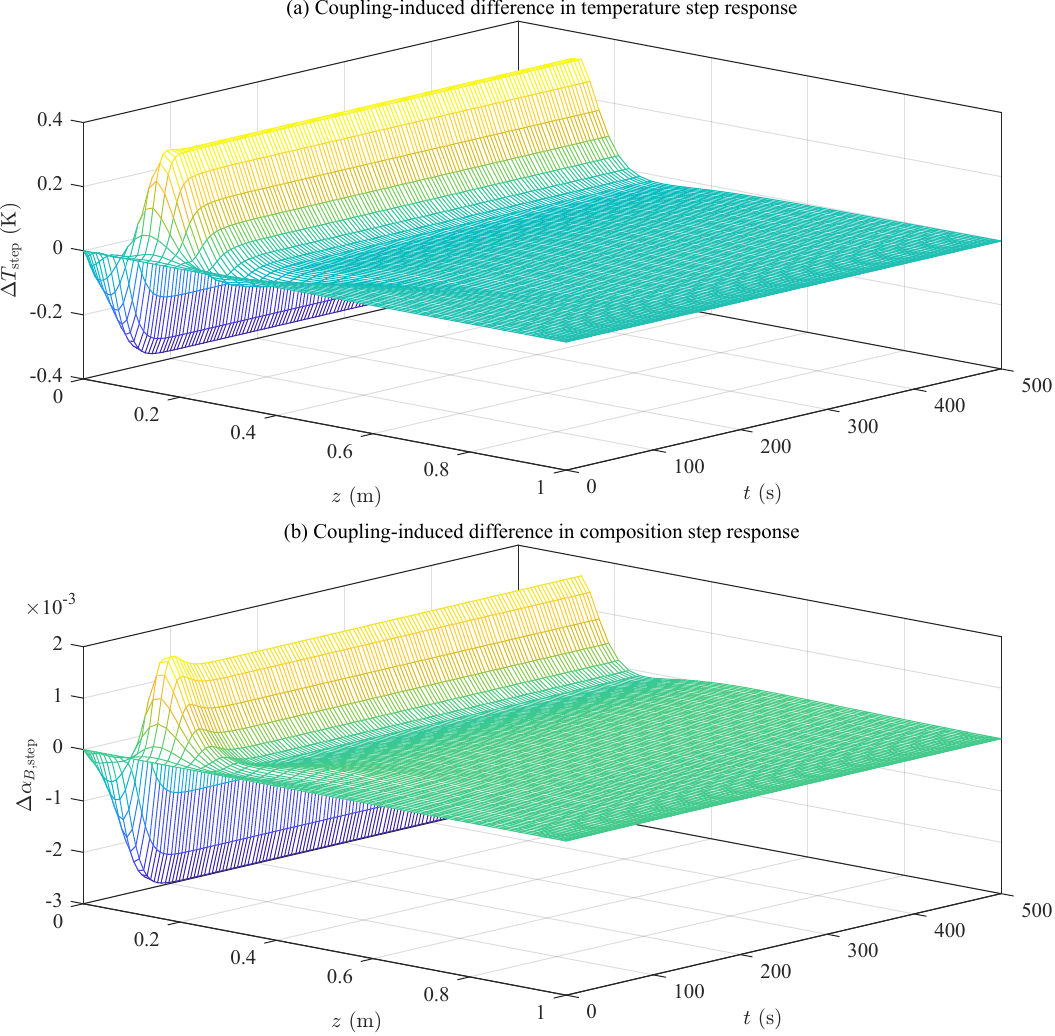}
    \caption{Effect of the non-diagonal Soret--Dufour coupling on the
    distributed jacket-step response:
    (a) $\Delta T_{\mathrm{step}}(z,t)$ and
    (b) $\Delta\alpha_{B,\mathrm{step}}(z,t)$,
    defined by \eqref{eq:coupling_step_difference}.
    The comparison uses $\chi=0.10$ and the diagonal case $\chi=0$,
    each initialized from its own $T_j=350~\mathrm{K}$ steady state.}
    \label{fig:soret_dufour_difference}
\end{figure}

The off-diagonal coupling produces a measurable modification of the
distributed transient, with
\[
    \max_{z,t}|\Delta T_{\mathrm{step}}|
    \simeq3.47\times10^{-1}~\mathrm{K},
    \qquad
    \max_{z,t}|\Delta\alpha_{B,\mathrm{step}}|
    \simeq2.76\times10^{-3}.
\]
The strongest differences occur in the upstream development region, where
both $|\partial_zT|$ and $|\partial_zA_r|$ are largest. This follows
directly from the coupled constitutive law: the Soret-type contribution to
$J$ is driven by the temperature gradient, whereas the Dufour-type
contribution to $q$ is driven by the affinity gradient. As these gradients
decrease downstream, the cross-transport contribution also weakens.

\begin{figure}[t]
    \centering
    \includegraphics[width=\columnwidth]
    {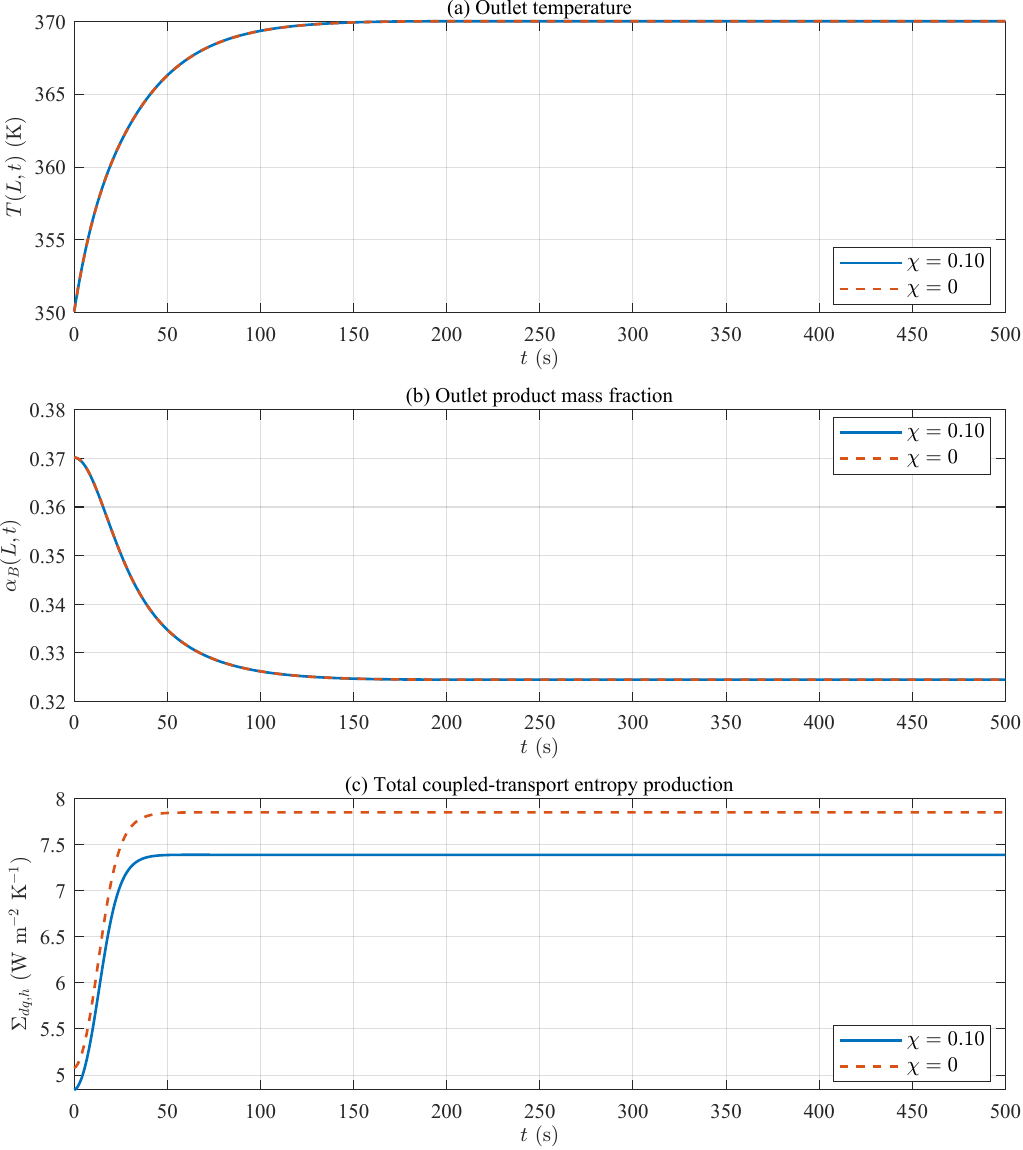}
    \caption{Comparison of the non-diagonal ($\chi=0.10$) and diagonal
    ($\chi=0$) transport models:
    (a) outlet temperature;
    (b) outlet product mass fraction; and
    (c) total diffusion--heat entropy production
    $\Sigma_{dq,h}(t)$.}
    \label{fig:soret_dufour_outlet_entropy}
\end{figure}

Although the internal distributed fields are measurably modified, the
outlet trajectories in
Fig.~\ref{fig:soret_dufour_outlet_entropy}(a)--(b) remain close. Thus,
under the present operating conditions, reciprocal diffusion--heat
coupling modifies the internal transport pathway more strongly than the
outlet observables.

For $\chi=0$, the total diffusion--heat entropy production reduces to
\begin{equation}
    \Sigma_{dq,h}^{\,\chi=0}
    =
    \Sigma_{d,h}+\Sigma_{q,h},
    \label{eq:diagonal_transport_entropy_relation}
\end{equation}
whereas for $\chi=0.10$ it is evaluated from
\eqref{eq:numerical_coupled_entropy}. Over the $500~\mathrm{s}$ interval,
including the respective steady baseline contributions,
\[
    \int_0^{500}
    \Sigma_{dq,h}^{\,\chi=0.10}(t)\,dt
    \simeq
    3.654\times10^3~
    \mathrm{J\,m^{-2}K^{-1}},
\]
compared with
\[
    \int_0^{500}
    \Sigma_{dq,h}^{\,\chi=0}(t)\,dt
    \simeq
    3.880\times10^3~
    \mathrm{J\,m^{-2}K^{-1}}.
\]
For this operating condition, the non-diagonal law therefore gives a lower
time-integrated diffusion--heat entropy production than the diagonal model.
This is not a general thermodynamic consequence of Soret--Dufour coupling:
positive semidefiniteness constrains the complete Onsager quadratic form to
be non-negative, whereas the cross contribution may have either sign.

\section{Conclusion}

A standard boundary-controlled irreversible port--Hamiltonian formulation
has been developed for a reacting tubular reactor with distributed
convection, chemical reaction, mass diffusion, and heat conduction. By
using an appropriate Eulerian thermodynamic coordinate representation, the
reactor dynamics are embedded in the standard BC-IPHS framework, with
thermodynamically admissible irreversible mechanisms and a
power-preserving boundary ports. A conservative-coordinate PFEM
discretization is further constructed, preserving the semi-discrete
Hamiltonian balance and admitting a discrete entropy balance with
non-negative internal entropy production.

The numerical study confirms the consistency of the spatial discretization
for a dimensional open-reactor benchmark and illustrates the distributed
reactor response in the intended low-Mach regime. The standard BC-IPHS
result applies to the diagonal diffusion--heat transport model, while a
separate Soret--Dufour comparison is included only to assess the effect of
non-diagonal reciprocal transport coupling. The resulting framework
provides a basis for future interconnection and boundary-control
developments for irreversible distributed reacting systems.

\appendix

\bibliography{biblio}

\begin{thebibliography}{10}
\expandafter\ifx\csname url\endcsname\relax
  \def\url#1{\texttt{#1}}\fi
\expandafter\ifx\csname urlprefix\endcsname\relax\def\urlprefix{URL }\fi
\expandafter\ifx\csname href\endcsname\relax
  \def\href#1#2{#2} \def\path#1{#1}\fi

\bibitem{maschke1992}
B.~M. Maschke, A.~J. van~der Schaft, Port controlled hamiltonian systems:
  Modelling origins and system theoretic properties, in: Proceedings of the
  IFAC Symposium on Nonlinear Control Systems Design (NOLCOS), 1992, pp.
  359--365.

\bibitem{vanderschaft2002}
A.~J. van~der Schaft, B.~M. Maschke, Hamiltonian formulation of
  distributed-parameter systems with boundary energy flow, Journal of Geometry
  and Physics 42~(1--2) (2002) 166--194.

\bibitem{legorrec2005}
Y.~Le~Gorrec, H.~Zwart, B.~Maschke, Dirac structures and boundary control
  systems associated with skew-symmetric differential operators, SIAM Journal
  on Control and Optimization 44~(5) (2005) 1864--1892.

\bibitem{maschke2013}
B.~Maschke, A.~J. van~der Schaft, P.~Crouch, Energy-based modeling and control
  of physical systems, in: Mathematical Control Theory I, Vol. 461 of Lecture
  Notes in Control and Information Sciences, Springer, Berlin, Germany, 2013,
  pp. 291--354.

\bibitem{Rashad2020}
R.~Rashad, F.~Califano, A.~J. van~der Schaft, S.~Stramigioli, Twenty years of
  distributed port-hamiltonian systems: a literature review, IMA Journal of
  Mathematical Control and Information 37~(4) (2020) 1400--1422.
\newblock \href {https://doi.org/10.1093/imamci/dnaa018}
  {\path{doi:10.1093/imamci/dnaa018}}.

\bibitem{ramirez2012}
H.~Ramirez, B.~Maschke, D.~Sbarbaro, Irreversible port hamiltonian systems,
  IFAC Proceedings Volumes 45~(19) (2012) 13--18.

\bibitem{Ramirez2013}
H.~Ramirez, B.~Maschke, D.~Sbarbaro, Irreversible port-hamiltonian systems: A
  general formulation of irreversible processes with application to the {CSTR},
  Chemical Engineering Science 89 (2013) 223--234.
\newblock \href {https://doi.org/10.1016/j.ces.2012.12.002}
  {\path{doi:10.1016/j.ces.2012.12.002}}.

\bibitem{ramirez2022overview}
H.~Ramirez, Y.~Le~Gorrec, An overview on irreversible port-hamiltonian systems,
  Entropy 24~(10) (2022) 1478.

\bibitem{ramirez2022bciphs}
H.~Ramirez, Y.~Le~Gorrec, B.~Maschke, Boundary controlled irreversible
  port-hamiltonian systems, Chemical Engineering Science 248 (2022) 117107.

\bibitem{zhou2017distributed}
W.~Zhou, B.~Hamroun, F.~Couenne, Y.~Le~Gorrec, Distributed port-hamiltonian
  modelling for irreversible processes, Mathematical and Computer Modelling of
  Dynamical Systems 23~(1) (2017) 3--22.

\bibitem{zhou2021thermodynamic}
W.~Zhou, B.~Hamroun, Y.~Le~Gorrec, F.~Couenne, A thermodynamic approach to the
  stabilization of tubular reactors, Journal of Process Control 108 (2021)
  98--111.

\bibitem{AltmannSchulze2017}
R.~Altmann, P.~Schulze, A port-hamiltonian formulation of the navier--stokes
  equations for reactive flows, Systems \& Control Letters 100 (2017) 51--55.
\newblock \href {https://doi.org/10.1016/j.sysconle.2016.12.005}
  {\path{doi:10.1016/j.sysconle.2016.12.005}}.

\bibitem{Ouardi2026}
A.~Ouardi, A.~Sarkar, H.~Ramirez, Y.~Le~Gorrec, Irreversible port-hamiltonian
  formulations for 1-dimensional fluid systems (2026).
\newblock \href {http://arxiv.org/abs/2603.11260} {\path{arXiv:2603.11260}}.

\bibitem{RamirezLeGorrec2024}
H.~Ramirez, Y.~Le~Gorrec, Interconnection of irreversible port hamiltonian
  systems, Automatica 170 (2024) 111846.
\newblock \href {https://doi.org/10.1016/j.automatica.2024.111846}
  {\path{doi:10.1016/j.automatica.2024.111846}}.

\bibitem{Kotyczka2016}
P.~Kotyczka, Finite volume structure-preserving discretization of 1d
  distributed-parameter port-hamiltonian systems, IFAC-PapersOnLine 49~(8)
  (2016) 298--303.
\newblock \href {https://doi.org/10.1016/j.ifacol.2016.07.457}
  {\path{doi:10.1016/j.ifacol.2016.07.457}}.

\bibitem{KotyczkaMaschkeLefevre2018}
P.~Kotyczka, B.~Maschke, L.~Lef{\`e}vre, Weak form of stokes--dirac structures
  and geometric discretization of port-hamiltonian systems, Journal of
  Computational Physics 361 (2018) 442--476.
\newblock \href {https://doi.org/10.1016/j.jcp.2018.02.006}
  {\path{doi:10.1016/j.jcp.2018.02.006}}.

\bibitem{CardosoRibeiro2018}
F.~L. Cardoso-Ribeiro, D.~Matignon, L.~Lef{\`e}vre, A structure-preserving
  partitioned finite element method for the 2d wave equation, IFAC-PapersOnLine
  51~(3) (2018) 119--124.
\newblock \href {https://doi.org/10.1016/j.ifacol.2018.06.033}
  {\path{doi:10.1016/j.ifacol.2018.06.033}}.

\bibitem{CardosoRibeiro2021}
F.~L. Cardoso-Ribeiro, D.~Matignon, L.~Lef{\`e}vre, A partitioned finite
  element method for power-preserving discretization of open systems of
  conservation laws, IMA Journal of Mathematical Control and Information 38~(2)
  (2021) 493--533.
\newblock \href {https://doi.org/10.1093/imamci/dnaa038}
  {\path{doi:10.1093/imamci/dnaa038}}.

\bibitem{Serhani2019}
A.~Serhani, D.~Matignon, G.~Haine, A partitioned finite element method for the
  structure-preserving discretization of damped infinite-dimensional
  port-hamiltonian systems with boundary control, in: Geometric Science of
  Information, 2019, pp. 549--558.
\newblock \href {https://doi.org/10.1007/978-3-030-26980-7_57}
  {\path{doi:10.1007/978-3-030-26980-7_57}}.

\bibitem{degroot1962}
S.~R. de~Groot, P.~Mazur, Non-Equilibrium Thermodynamics, North-Holland,
  Amsterdam, The Netherlands, 1962.

\bibitem{vanbrunt2022}
J.~Van~Brunt, E.~G. Piedrahita, B.~E. Larson, D.~A. Miller,
  Onsager--stefan--maxwell formulation of non-isothermal multicomponent
  diffusion with soret and dufour effects, AIChE Journal 68~(3) (2022) e17625.

\bibitem{ToledoZucco2024}
J.-P. Toledo-Zucco, D.~Matignon, C.~Poussot-Vassal, Y.~Le~Gorrec,
  Structure-preserving discretization and model order reduction of
  boundary-controlled 1d port-hamiltonian systems, Systems \& Control Letters
  194 (2024) 105947.
\newblock \href {https://doi.org/10.1016/j.sysconle.2024.105947}
  {\path{doi:10.1016/j.sysconle.2024.105947}}.

\bibitem{KotyczkaLefevre2019}
P.~Kotyczka, L.~Lef{\`e}vre, Discrete-time port-hamiltonian systems: A
  definition based on symplectic integration, Systems \& Control Letters 133
  (2019) 104530.
\newblock \href {https://doi.org/10.1016/j.sysconle.2019.104530}
  {\path{doi:10.1016/j.sysconle.2019.104530}}.

\end{thebibliography}

\end{document}